\documentclass[12pt]{amsart}
\usepackage{amssymb}

\begin{document}

\numberwithin{equation}{section}

\newtheorem{theorem}{Theorem}[section]
\newtheorem{proposition}[theorem]{Proposition}
\newtheorem{conjecture}[theorem]{Conjecture}
\def\theconjecture{\unskip}
\newtheorem{corollary}[theorem]{Corollary}
\newtheorem{lemma}[theorem]{Lemma}
\newtheorem{observation}[theorem]{Observation}
\theoremstyle{definition}
\newtheorem{definition}{Definition}
\def\thedefinition{\unskip}
\newtheorem{remark}{Remarks}
\def\theremark{\unskip}
\newtheorem{question}{Question}
\def\thequestion{\unskip}
\newtheorem{example}{Example}
\def\theexample{\unskip}
\newtheorem{problem}{Problem}

\def\intprod{\mathbin{\lr54}}
\def\reals{{\mathbb R}}
\def\integers{{\mathbb Z}}
\def\complex{{\mathbb C}\/}
\def\distance{\operatorname{distance}\,}
\def\ZZ{ {\mathbb Z} }
\def\e{\varepsilon}
\def\p{\partial}
\def\rp{{ ^{-1} }}
\def\Re{\operatorname{Re\,} }
\def\Im{\operatorname{Im\,} }

\def\intprod{\mathbin{\lr54}}
\def\reals{{\mathbb R}}
\def\integers{{\mathbb Z}}
\def\complex{{\mathbb C}\/}
\def\distance{\operatorname{distance}\,}
\def\ZZ{ {\mathbb Z} }
\def\e{\varepsilon}
\def\p{\partial}
\def\rp{{ ^{-1} }}
\def\Re{\operatorname{Re\,} }
\def\Im{\operatorname{Im\,} }

\def\scriptx{{\mathcal X}}
\def\scripti{{\mathcal I}}
\def\scripth{{\mathcal H}}
\def\scriptm{{\mathcal M}}
\def\scripte{{\mathcal E}}
\def\scriptt{{\mathcal T}}
\def\scriptb{{\mathcal B}}
\def\frakg{{\mathfrak g}}
\def\frakG{{\mathfrak G}}

\def\ov{\overline}

\author{Michael Christ}
\address{
        Michael Christ\\
        Department of Mathematics\\
        University of California \\
        Berkeley, CA 94720-3840, USA}
\email{mchrist@math.berkeley.edu}
\thanks{
The first author was supported by NSF grant DMS-9970660,
and by the Miller Institute for Basic Research in Science,
while he was on appointment as a Miller Research Professor}

\author{Alexander Kiselev}
\address{ Alexander Kiselev\\
Department of Mathematics\\
University of Chicago\\
Chicago, Ill. 60637}
\email{kiselev@math.uchicago.edu}
\thanks{The second author was supported by NSF grant DMS-0102554
and by an Alfred P. Sloan Research Fellowship}

\title[Scattering and wave operators]
{Scattering and wave operators for
one-dimensional Schr\"odinger operators with slowly
decaying nonsmooth potentials}

\begin{abstract}
We prove existence of modified wave operators for one-dimensional Schr\"odinger equations
with potential in $L^p(\reals)$, $p<2$. If in addition the potential is 
conditionally integrable, then the usual   M\"oller wave operators exist. 
We also prove asymptotic completeness of these wave operators for 
some classes of random potentials, and for almost every boundary condition
for any given potential.
\end{abstract}
\date{\today}
\maketitle


\section{Introduction}

Let $H_V$ be a one-dimensional Schr\"odinger operator defined by
\begin{equation}\label{start}
H_V =  -\frac{d^2}{dx^2} + V(x). 
\end{equation}
Let us discuss the case where the operator is defined on a half-axis, with some 
self-adjoint boundary condition at zero. 
We are interested in potentials decaying at infinity, for which we may expect that 
asymptotically as time tends to infinity, motion of the associated perturbed quantum system
resembles the free evolution.
What is the critical rate of
decay of the potential for which the dynamics remains close to free for large times?
The mathematical framework for studying these questions is provided
 by scattering theory. 
Recall that the wave operators associated with $H_V$ and $H_0$ are defined by 
\begin{equation}\label{wo}
\Omega_\pm f = \lim_{t \rightarrow \mp \infty} e^{itH_V} e^{-itH_0} f,
\end{equation}
where the limit is understood in the strong $L^2$ sense. In particular,
 existence of the 
wave operators implies that the absolutely continuous spectrum of $H_V$
 fills all of $\reals^+$ (see, e.g. \cite{RS3})
and moreover provides significant information on large time dynamics 
$e^{-itH_V}.$ 

The wave operators will be called asymptotically complete if the range of $\Omega_\pm$ 
coincides with the orthogonal complement of the subspace spanned by eigenvectors of
the operator $H_V$. An alternative equivalent characterization is that the range of 
the wave operators
is equal to the absolutely continuous subspace $\scripth_{\text{ac}}(H_V)$ of the operator
$H_V$, and that the singular continuous spectrum $\sigma_{{\rm sc}}(H_V)$ is empty.  
We note that the intended intuitive meaning of asymptotic completeness 
is that the dynamics 
of the perturbed operator can be divided into two well-understood parts: 
scattering states traveling 
to infinity in a way similar to the free evolution, and bound states which
 remain confined in a certain sense
for all times.  
We postpone more detailed discussion of
the notion of asymptotic completeness to Section \ref{section:completeness}.

A well-known result, going back to the very first years of the rigorous
 scattering theory, says that 
if $V \in L^1,$ then the wave operators exist and are asymptotically complete. 
Moreover, in this case the spectrum on the positive 
semi-axis is purely absolutely continuous, and there can be only 
discrete spectrum below zero, possibly accumulating at zero.
There has been much work extending the existence of wave operators to wider
classes of potentials,
with some additional conditions on derivatives,
or for potentials of certain particular forms.
Generally, for more slowly decaying 
potentials one needs to modify the free dynamics in the definition of wave
operators in order for the limits to exist.
The first work of this type was that of Dollard \cite{Dol}, who studied the 
Coulomb potential.  Further developments used
computation of asymptotic classical trajectories by means of a Hamilton-Jacobi equation
to build an appropriate phase correction to the free dynamics, which was used to prove
existence of modified wave operators (in any dimension).  
See, for instance,  Buslaev and Matveev \cite{BM}, Ahsholm and Kato \cite{AK}; the
 strongest results are contained in 
H\"ormander \cite{hormanderwaveoperators}. For example, 
H\"ormander \cite{hormanderwaveoperators} gives 
existence of wave operators if $|V(x)| \leq C(1+|x|)^{-1/2-\epsilon}$ and
$|D^{\alpha}V(x)| \leq C(1+|x|)^{-3/2-\epsilon}$ for any $|\alpha|=1.$ 
One needs more conditions on derivatives to compensate for slower decay.

Another series of works (see \cite{BA,BAD,D1,W1} for further references)
 studied oscillating potentials of 
Wigner-von Neumann type, for example $V(x) = \sin (cx^{\alpha})/x^{\beta}.$
Interest in this class
of potentials was in part fueled by the Wigner-von Neumann example \cite{WvN},
 where $V(x) \sim c\sin (2x)/x$ at infinity,
which leads to a positive eigenvalue embedded in the absolutely continuous
 spectrum. 
For a wide class of potentials of this type (including the original example of \cite{WvN}),
existence and asymptotic completeness of (sometimes modified) wave operators has been shown.

In the opposite direction, Pearson \cite{Pe} constructed examples of
potentials in $\cap_{p>2} L^p(\reals)$ for which the spectrum 
is purely singular, and hence wave operators do not exist. Kotani and
 Ushiroya \cite{KU} also provided 
power decaying examples where the spectrum is purely singular (pure point) for the
 rate of decay $x^{-\alpha},$ $\alpha<1/2.$ 

In the last few years, there has been a series of works 
studying  existence of absolutely continuous spectrum for slowly decaying 
potentials in full generality, with no additional conditions on derivatives
or specific representation. In \cite{christkiselevpowerdecay,Re1} it was shown 
that the absolutely continuous spectrum is preserved if $|V(x)| \leq 
C(1+x)^{-\alpha},$ $\alpha > 1/2,$ and moreover the generalized 
eigenfunctions have WKB-type plane wave asymptotics. 
Recently, we improved the result of \cite{christkiselevpowerdecay} to treat potentials 
in $L^p,$ $p<2$ \cite{christkiselevdecaying}. 
The sharpest result on the stability of the absolutely continuous spectrum 
is due to Deift and Killip \cite{DK}, who showed it for $V \in L^2$. 
This result is optimal in the $L^p$ scale, as Pearson's examples show. 
The method of \cite{DK} is quite different from \cite{christkiselevpowerdecay,Re1} and 
yields much less information concerning the nature of the generalized eigenfunctions. 

Although the absolute continuity of the spectrum is an important
characterization of an operator with direct consequences for the 
dynamical behavior of the particle, wave operators provide a much finer
description of long-time dynamics. 
The purpose of this paper is to use the additional information contained
in the asymptotic behavior of generalized eigenfunctions \cite{christkiselevdecaying} to prove 
the existence of (modified) wave operators for potentials $V \in L^p,$
$p<2.$
Let
\begin{equation}\label{phasec}
W(\lambda,t) = -(2\lambda)\rp\int_0^{2\lambda t} V(s)\,ds.
\end{equation}
Define $e^{-iH_0t\pm iW(H_0,\mp t)}$ to be the Fourier multiplier
operator on $L^2(\reals^+)$ that maps
$\int_0^\infty F(\lambda)\sin(\lambda x)\,d\lambda$
to $\int_0^\infty e^{-i\lambda^2t \pm iW(\lambda^2,\mp t)}
F(\lambda)\sin(\lambda x)\,d\lambda$,
for all $F\in L^2(\reals^+,d\lambda)$.
 Define
\begin{equation} \label{modifiedwaveoperatorsdefn}
\Omega^m_\pm f
= \lim_{t\to\mp\infty} e^{itH_V}e^{-it H_0 \pm iW(H_0,\mp t)}f
\end{equation}
for all $f\in L^2(\reals^+)$.
Among our main conclusions are the following two theorems.
\begin{theorem}  \label{thm:waveoperators}
Let $V$ be a potential in the class $L^1+L^p(\reals^+)$ 
for some $1<p<2$,
and let $H_V$ be the associated Schr\"odinger operator on
$L^2(\reals^+)$ with any self-adjoint boundary condition at $0$.
Then for every $f\in L^2(\reals^+)$, the limits in 
\eqref{modifiedwaveoperatorsdefn} exist in $L^2(\reals^+)$ norm,
as $t\to\mp\infty$.
The modified wave
operators $\Omega^m_\pm$ thus defined are both unitary bijections
from $\scripth=L^2(\reals^+)$ to $\scripth_{\text{ac}}(H_V)$.
\end{theorem}

\begin{theorem}\label{mollerwo}
In addition to the hypotheses of the preceding theorem,
suppose that the improper integral $\int_0^\infty V(s)\,ds$ exists. 
Then for every $f\in L^2(\reals^+)$,
the limits in \eqref{wo} exist in $L^2(\reals^+)$ norm,
as $t\to\mp\infty$.
The wave operators $\Omega_\pm$ thus defined are both unitary bijections
from $\scripth=L^2(\reals^+)$ to $\scripth_{\text{ac}}$.
\end{theorem}
By the improper integral we mean of course $\lim_{N\to\infty}\int_0^N V$; we are not
assuming that $V\in L^1$.  
We also prove analogous results on wave operators 
for the whole line case and 
for some Dirac operators, see Section \ref{section:modifiedwave}.

Another way to put the conclusion is that for each $f\in L^2(\reals^+)$
there exist functions $F_\pm\in \scripth_{\text{ac}}$ such that
$\|e^{itH_V}F_\pm-e^{itH_0}f\|_{L^2(\reals^+)}\to 0$
as $t\to\pm\infty$, and each map $f\mapsto F_\pm$ is unitary
and surjective onto $\scripth_{\text{ac}}$.
Thus there exists a full family of scattering states, and any state that is 
asymptotically free at $t=\mp\infty$ is likewise free at $t = \pm\infty$.

The main idea of the proof is to use generalized eigenfunctions 
to construct the spectral decomposition of the operator $H_V$, 
and in particular to 
derive an explicit expression for the evolution group 
of the perturbed operator.
Generalized eigenfunctions play a role parallel to the Fourier transform 
in the free case. The existence of wave operators is then proven 
by direct analysis, comparing the perturbed evolution with modified 
free dynamics. While the analysis is generally based on results of
\cite{christkiselevpowerdecay,christkiselevdecaying}, we need 
several essential improvements in the estimates
to fulfill this plan. In Sections 
\ref{section:numerical}, \ref{section:multilinear} 
we extend the 
analysis of multilinear expressions encountered in the series for 
generalized eigenfunctions to cover the situations arising in applications
to wave operators.  
In Section \ref{section:complex} we consider solutions of $H_V u = zu$ 
for complex energies $z$, and establish certain uniform bounds and asymptotics as $z\to\reals^+$.
In Section \ref{section:limiting absorption} we prove a limiting absorption principle, 
which allows us to write an explicit formula for the evolution. 
Sections \ref{section:heart}, \ref{section:phasereduction} and \ref{section:modifiedwave} 
contain long-time asymptotic analysis and the rest 
of the proof of the existence of wave operators. 
The key step is Lemma~\ref{truekeylemma}, where the full strength of the multilinear analysis
is required to justify discarding all terms with which it is concerned.
In Section \ref{section:completeness} we discuss the issue of asymptotic completeness for 
generic potentials within certain classes.

Although all theorems, and many of the intermediate results, 
are stated for potentials in $L^1+L^p$ for some $1<p<2$,
for simplicity we first give the proofs under the simplifying assumption that
$V\in L^p$.
In the last section we review the machinery \cite{christkiselevdecaying}
needed to extend the analysis from $L^p$ to $L^1+L^p$, and to the more general 
amalgamated class $\ell^p(L^1)$.

We note that the existence of (modified) wave operators
for general long-range potentials in higher dimensions remains
an open problem. For recent progress in this direction, see
\cite{rodnianskischlag}.

Some of the results of this paper were announced in \cite{ipamnotes}.

\section{A numerical bound for iterated multiple integrals}
\label{section:numerical}

Let $\{f_i\}$ be a collection of integrable functions
from $\reals$ to $\complex$.
Consider multilinear expressions 
\begin{equation}
M_n(f_1,\dots,f_n)
= \iint_{x_1\le\cdots\le x_n}
\prod_{i=1}^n f_i(x_i)dx_i
\end{equation}
and their maximal variants
\begin{equation}
M_n^*(f_1,\dots,f_n)
= \sup_{y}
\Big| \iint_{x_1\le\cdots\le x_n\le y}
\prod_{i=1}^n f_i(x_i)dx_i \Big| . 
\end{equation}
The purpose of this section is to establish upper bounds
for $M_n,M_n^*$, in terms of certain auxiliary functions
$g_\delta$ of the functions $f_i$, with particular attention
to the dependence on $n$ as $n\to\infty$. These bounds will 
play an essential part in our analysis of asymptotics
of wave groups.

\begin{definition}
A martingale structure $\{E^m_j\}$ on an interval $I\subset\reals$
is a collection of subintervals $E^m_j\subset I$,
indexed by $m\in\{0,1,2,\dots\}$ and $j\in\{1,2,\dots,2^m\}$,
possessing the following two properties.
(i) Except for endpoints, $\{E^m_j: 1\le j\le 2^m\}$ is
a partition of $I$, for each $m$.
(ii) $E^m_j = E^{m+1}_{2j-1}\cup E^{m+1}_{2j}$ for all $m,j$.
\end{definition}
To any $f\in L^1$, any $\delta\in\reals$,
and any martingale structure, we associate 
\begin{equation}  \label{eq:defnofgdelta}
g_\delta(f) = 
\sum_{m=1}^\infty 2^{\delta m}
(\sum_{j=1}^{2^m} |{\textstyle\int}_{E^m_j}f|^2)^{1/2} .
\end{equation}
More generally, define
\begin{equation}
g_\delta(\{f_k\}) = 
\sum_{m=1}^\infty 2^{\delta m}
(\sum_{j=1}^{2^m} \sup_k |{\textstyle \int}_{E^m_j}f_k|^2)^{1/2} .
\end{equation}


\begin{proposition}  \label{prop:numericalbound}
There exists $C<\infty$ such that
for any martingale structure $\{E^m_j\}$,
any $\delta\ge 0$,
any $f_1,\dots,f_n\in L^1(\reals)$,
and any $n\ge 2$,
\begin{equation}
|M_n(f_1,\dots,f_n)|
\le \frac{C^{n+1}}{\sqrt{n!}}
g_{-\delta}(f_1)\cdot g_\delta(\{f_k:k\ge 2\})^{n-1}.
\end{equation}
Moreover for any $\delta'>\delta \ge 0$,
there exists $C<\infty$ such that for all $\{f_i\}$ and all $n\ge 2$,
\begin{equation}
|M_n^*(f_1,\dots,f_n)|
\le \frac{C^{n+1}}{\sqrt{n!}}
g_{-\delta}(f_1)\cdot g_{\delta'}(\{f_k:k\ge 2\})^{n-1}.
\end{equation}
\end{proposition}

In previous work \cite{christkiselevfiltrations}
we proved the simpler analogue 
with $\delta=0$ and with the bound (for $M_n$)
$C^{n+1}(n!)^{-1/2}g(\{f_k: k\ge 1\})^n$, 
and applied it to the analysis of generalized eigenfunctions,
which can be expanded as sums over $n$ of such
iterated multiple integrals, where $f_k$ is
essentially $\exp(\pm 2i\lambda x)V(x)$, $V$ is the potential,
and $\lambda^2$ is a spectral parameter;
$g(f)$ is then a function of $\lambda$.
In the present work, we need a refinement in which
$f_1$ is essentially $\exp(\pm i\lambda x)h(x)$,
and $h$ is an arbitrary $L^2$ function, unrelated to $V$.
The quantity $g(h)$ is not appropriately bounded for
$h\in L^2$, forcing the introduction of the mollifying
factors $2^{-\delta m}$ in its definition. This in turn
forces compensating factors of $2^{+\delta m}$,
leading to the above formulation.

\begin{proof}
It is proved in \cite{christkiselevfiltrations}
that there exist positive constants $b_n$ satisfying
$b_n\le C^{n+1}/\sqrt{n!}$
and $n^{1/2}b_{n+1}/b_n\to C$ as $n\to\infty$,
such that for all nonnegative real numbers $x,y$,
\begin{equation} \label{eq:modifiedbinomial}
b_n y^n + \sum_{i=2}^{n-2} b_i b_{n-i}x^iy^{n-i} + b_n x^n
\le b_n (x^2+y^2)^{n/2}.
\end{equation}
It is also shown in \cite{christkiselevfiltrations} that
\begin{equation} \label{eq:onefunctionversion}
\begin{split}
|M_n(f,\dots,f)| &\le b_n g_0(f)^n
\\
M_n^*(f,\dots,f) &\le C_\e^n b_n g_\e(f)^n
\end{split}
\end{equation}
for every $\e>0$. Moreover,
for distinct functions $f_i$,
\begin{equation} \label{eq:modifiedonefunctionversion}
\begin{split}
|M_n(f_1,\dots,f_n)| &\le b_n g_0(\{f_i\})^n
\\
M_n^*(f_1,\dots,f_n)|&\le C_\e^n b_n g_\e(\{f_i\})^n.
\end{split}
\end{equation}
Although this bound is not explicitly formulated in 
\cite{christkiselevfiltrations}, it follows directly from exactly  
the argument given there.

For $n\ge 2$ define $\tilde b_n = R^{n} b_{n-2}$,
where $R$ is a sufficiently large positive constant, to be
determined later in the proof.
In order to simplify notation, we will prove
the result in the special case $f_2=f_3=\cdots =f_n=f$,
and will write $f_1=\tilde f$.
The proof will be by induction on $n$. First we will treat only
the case where $n\ge 6$, assuming the result for 
all $n\le 5$, and at the end will discuss the modification
for small $n$.

We write 
$f^m_j = \chi_{E^m_j}\cdot f$ and $\tilde f^m_j 
= \chi_{E^m_j}\cdot \tilde f$,
where $\chi_{E}$ denotes always the characteristic function of a set $E$.
\begin{lemma}
If $R$ is chosen to be sufficiently large then for all $n\ge 2$
and any $\delta\ge 0$,
\begin{equation}
|M_n(\tilde f,f,\dots,f)|\le \tilde b_n g_{-\delta}(\tilde f)
\cdot g_\delta(f)^{n-1}.
\end{equation}
\end{lemma}

\begin{proof}
By inequality (4.6) of \cite{christkiselevfiltrations},
\begin{equation}\begin{split} \label{eq:inductivestep}
|M_n(\tilde f,f,\dots,f)|
\le 
&|M_n(\tilde f^1_1,f^1_1,\dots,f^1_1)|
\\
& + |M_{n-1}(\tilde f^1_1,f^1_1,\dots,f^1_1)|\cdot 
|{\textstyle\int} f^1_2|
\\
& + \sum_{i=2}^{n-2} |M_{n-i}(\tilde f^1_1,f^1_1,\dots,f^1_1)|
\cdot |M_i(f^1_2,\dots,f^1_2)|
\\
& + |{\textstyle\int}\tilde f^1_1|\cdot M_{n-1}(f^1_2,\dots,f^1_2)|
+ |M_n(\tilde f^1_2,f^1_2,\dots,f^1_2)|
\\
\le
& |M_n(2^{-\delta}\tilde f^1_1,2^\delta f^1_1,\dots,2^\delta f^1_1)|
\\
& + |M_{n-1}(2^{-\delta} \tilde f^1_1,
2^\delta f^1_1,\dots,2^\delta f^1_1)|
\cdot |{\textstyle\int} 2^\delta f^1_2|
\\
& + \sum_{i=2}^{n-2} |M_{n-i}( 2^{-\delta}\tilde f^1_1,2^\delta f^1_1,\dots,
2^\delta f^1_1)|
\cdot |M_i(2^\delta f^1_2,\dots,2^\delta f^1_2)|
\\
& + |{\textstyle\int} 2^{-\delta}\tilde f^1_1|
\cdot |M_{n-1}(2^\delta f^1_2,\dots,2^\delta f^1_2)|
\\
& + |M_n(2^{-\delta}\tilde f^1_2,2^\delta f^1_2,\dots,2^\delta f^1_2)| .
\end{split}\end{equation}
We have assumed that $n\ge 2$ and $\delta\ge 0$ to ensure that
at least one factor of $2^{\delta}$ offsets the factor of $2^{-\delta}$.

The first and last terms in the preceding bound involve $M_n$ itself,
but the former involves only the restrictions of 
$\tilde f,f$ to $E^1_1$, while the latter involves only their restrictions
to $E^1_2$; thus these expressions are in a sense simpler than
the original expression $M_n$. We will therefore use as 
part of our induction hypothesis the desired inequalities for 
$M_n(\tilde f^1_1,f^1_1,\cdots,f^1_1)$ and for
$M_n(\tilde f^1_2,f^1_2,\cdots,f^1_2)$.
For the justification of this method of reasoning see
the two paragraphs immediately following inequality (4.12) of
\cite{christkiselevfiltrations}. 

The collection of all those sets $E^m_j$
with $m\ge 1$ and $j\le 2^{m-1}$
forms a martingale structure on $E^1_1$; 
however, when it is considered as such, the index $m$ should be 
replaced by $m-1$. Thus 
the induction hypothesis, for the first term on the right-hand
side of the preceding bound, 
may be stated as
\begin{multline}
|M_n(2^{-\delta} \tilde f^1_1,2^{\delta} f^1_1,\cdots, 2^\delta f^1_1)|
\\
\le \tilde b_n
\sum_{m=2}^\infty
2^{-m\delta}
(\
\sum_{j=1}^{2^{m-1}}
|{\textstyle\int} \tilde f^m_j|^2
)^{1/2}
\cdot
\Big[
\sum_{m=2}^\infty
2^{m\delta}
(\
\sum_{j=1}^{2^{m-1}}
|{\textstyle\int} f^m_j|^2
)^{1/2}\Big]^{n-1}\ .
\end{multline}
There is a corresponding bound for
$|M_n(2^{-\delta}\tilde f^1_2,2^\delta f^1_2,\cdots,2^\delta f^1_2)|$,
with the inner sum running instead over $2^{m-1}<j\le 2^m$.

To formulate the induction hypothesis for the general
term of the preceding expression, introduce
\begin{equation}
\tilde g^1_t 
= 
\sum_{m=2}^\infty 
2^{-\delta m}
\Big( \sum_{j: E^m_j\subset E^1_t}
|{\textstyle\int} \tilde f^m_j|^2
\Big)^{1/2}
\ \text{ and }\ 
g^1_t 
= 
\sum_{m=2}^\infty 
2^{\delta m}
\Big( \sum_{j: E^m_j\subset E^1_t}
|{\textstyle\int} f^m_j|^2
\Big)^{1/2}
\end{equation}
for $t=1,2$
and
\begin{equation} 
\tilde g^1  = ((\tilde g^1_1)^2 + (\tilde g^1_2)^2)^{1/2},
\qquad
g^1  = ((g^1_1)^2 + (g^1_2)^2)^{1/2} .
\end{equation}
Note that 
$(g^1)^2  + 2^{2\delta} |\int f^1_1|^2 + 2^{2\delta} |\int f^1_2|^2 
\leq g_\delta(f)^2$,
with a corresponding relation between 
$\tilde g^1,\,g_{-\delta}(\tilde f)$.

>From the induction hypothesis and \eqref{eq:modifiedonefunctionversion} 
we obtain 
\begin{equation}
\begin{split}
|M_n(\tilde f,f,\dots,f)|
\le
& \tilde b_n \tilde g^1_1(g^1_1)^{n-1}
+ \tilde b_{n-1}\tilde g^1_1(g^1_1)^{n-2} 
|{\textstyle\int} 2^\delta f^1_2| 
\\ 
& + \sum_{i=2}^{n-2} \tilde b_{n-i}b_i \tilde g^1_1(g^1_1)^{n-1-i} (g^1_2)^i
 + |{\textstyle\int} 2^{-\delta}\tilde f^1_1| \cdot b_{n-1}(g^1_2)^{n-1}
 + \tilde b_n \tilde g^1_2(g^1_2)^{n-1}.
\end{split}\end{equation}
Consider now the sum of the first term on the right, together with
all terms of the summation for which the index $i$ is either in $[2,n-4]$,
or equals $n-2$:
\begin{align*}
\tilde b_n \tilde g^1_1(g^1_1)^{n-1}
& + \sum_{i=2}^{n-4} \tilde b_{n-i}b_i \tilde g^1_1(g^1_1)^{n-1-i}
(g^1_2)^i
+ \tilde b_2 b_{n-2}\tilde g^1_1 g^1_1(g^1_2)^{n-2}
\\
& = R^n\tilde g^1_1 g^1_1
\Big(
b_{n-2}(g^1_1)^{n-2}
+ \sum_{i=2}^{n-4}R^{-i} b_{n-i-2}b_i (g^1_1)^{n-2-i}(g^1_2)^i
+ R^{-n+2}b_0 b_{n-2}(g^1_2)^{n-2}
\Big)
\\
& \le R^n\tilde g^1_1 g^1_1
\Big(
b_{n-2}(g^1_1)^{n-2}
+ \sum_{i=2}^{n-4}b_{n-i-2}b_i (g^1_1)^{n-2-i}(g^1_2)^i
+ b_0 b_{n-2}(g^1_2)^{n-2}
\Big)
\\
& \le R^n\tilde g^1_1 g^1_1
b_{n-2} 
\big(
(g^1_1)^{2(n-2)}+(g^1_2)^{2(n-2)}
\big)^{1/2}
\\
& = \tilde b_n \tilde g^1_1 g^1_1 (g^1)^{n-2} .
\end{align*}
To pass from the third line to the fourth we have invoked
\eqref{eq:modifiedbinomial}.

We now have
\begin{equation}  \label{intermediateboundforMn}
\begin{split}
|M_n(\tilde f,f,\dots,f)|
& \le \tilde b_n \tilde g^1_1 g^1_1 (g^1)^{n-2}
+ \tilde b_3 b_{n-3} \tilde g^1_1(g^1_1)^2(g^1_2)^{n-3}
+ \tilde b_n\tilde g^1_2(g^1_2)^{n-1}
\\
&+ \tilde b_{n-1}\tilde g^1_1(g^1_1)^{n-2}|{\textstyle\int} 2^\delta f^1_2|
+ |{\textstyle\int} 2^{-\delta}\tilde f^1_1| b_{n-1}(g^1_2)^{n-1} .
\end{split}\end{equation}
Set $\beta = \tilde b_3 b_{n-3}/\tilde b_n = R^{3-n}b_1b_{n-3}/b_{n-2}$ 
and note that $\beta\le CR^{3-n}n^{1/2}\le 1$, if $R$ is chosen
to be sufficiently large, since we are assuming $n\ge 6$.
Applying Cauchy-Schwarz to the sum of the first three terms on
the right-hand side of the preceding inequality, we obtain
\begin{multline}
\tilde b_n \tilde g^1_1 g^1_1 (g^1)^{n-2}
+ \tilde b_3 b_{n-3} \tilde g^1_1(g^1_1)^2(g^1_2)^{n-3}
+ \tilde b_n\tilde g^1_2(g^1_2)^{n-1}
\\
\le\tilde b_n\tilde g^1\cdot
\Big(
(g^1_1)^2(g^1)^{2n-4}
+ (g^1_2)^{2n-2}
+ 2\beta(g^1_1)^3(g^1_2)^{n-3}(g^1)^{n-2}
+ \beta^2(g^1_1)^4(g^1_2)^{2n-6}
\Big)^{1/2}.
\end{multline}

We claim that
\begin{equation}
(g^1_1)^2(g^1)^{2n-4}
+ (g^1_2)^{2n-2}
+ 2\beta(g^1_1)^3(g^1_2)^{n-3}(g^1)^{n-2}
+ \beta^2 (g^1_1)^4(g^1_2)^{2n-6}
\le (g^1)^{2n-2},
\end{equation}
provided that $n\ge 6$ and $\beta$ is sufficiently small.
Writing $x=g^1_1,y=g^1_2$, this is simply
\begin{equation*}
x^2(x^2+y^2)^{n-2}+2\beta x^3y^{n-3}(x^2+y^2)^{(n-2)/2}
+\beta^2 x^4 y^{2n-6}+y^{2n-2}
\le (x^2+y^2)^{n-1},
\end{equation*}
which is equivalent to
\begin{equation*}
2\beta x^3 y^{n-3}(x^2+y^2)^{(n-2)/2}
+\beta^2 x^4 y^{2n-6}+y^{2n-2}
\le y^2(x^2+y^2)^{n-2}.
\end{equation*}
By homogeneity we may assume that $x^2+y^2=1$,
so we wish to have
\begin{equation*}
2\beta x^3 y^{n-3}
+\beta^2 x^4 y^{2n-6}+y^{2n-2}
\le y^2.
\end{equation*}
This holds whenever $x^2+y^2=1$ and $x,y\ge 0$,
provided $\beta$ is sufficiently small,
provided that each of the exponents on $y$ on the left-hand side 
is strictly larger than $2$; this holds provided $n\ge 6$.

We have thus established that for $n\ge 6$,
\begin{multline}  \label{eq:lastMnstep1}
|M_n(\tilde f,f,\dots,f)|
\le\tilde b_n\tilde g^1(g^1)^{n-1}
\\
+ \tilde b_{n-1}\tilde g^1_1(g^1_1)^{n-2}|{\textstyle\int} 2^\delta f^1_2|
+ |{\textstyle\int} 2^{-\delta}\tilde f^1_1| b_{n-1}(g^1_2)^{n-1} .
\end{multline}
Representing the right-hand side as $\tilde b_n\tilde x x^{n-1}
+ \tilde b_{n-1}\tilde x x^{n-2}y + b_{n-1}\tilde y x^{n-1}$
where $x=g^1,y=|\int 2^\delta f^1_2|$
and $\tilde x = \tilde g^1,\tilde y 
= |{\textstyle\int} 2^{-\delta}\tilde f^1_1|$,
we seek to bound this right-hand side by
\begin{equation} \label{eq:lastMnstep2}
\tilde b_n (\tilde x+\tilde y)(x+y)^{n-1}
\ge \tilde b_n \tilde x x^{n-1} + \tilde b_n\tilde x(n-1)x^{n-2}y
+ \tilde b_n\tilde y x^{n-1}.
\end{equation}
Now\footnote{It is in order to be able to assert this last inequality
that we incorporate an $\ell^1$ with respect to $m$,
rather than an $\ell^2$ norm, in the definitions of the
functionals $g_\delta$.}
 $\tilde b_{n-1} \lesssim R n^{1/2}\tilde b_n\le(n-1)\tilde b_n$
for all sufficiently large $n$.
Likewise, $b_{n-1} \le Cn^{1/2}R\rp\tilde b_n
\le (n-1)\tilde b_n$.
Hence a term-by-term comparison completes the proof of the bound 
for $M_n$, for $n\ge 6$.
\end{proof}

Observe that the restriction $n\ge 6$ was only used above to
control the first line of the right-hand side of
\eqref{intermediateboundforMn}; once that term is majorized by
$\tilde b_n\tilde g^1(g^1)^{n-1}$, 
the reasoning of the final paragraph allows us to absorb the two
remaining special terms.

The cases $n=2,3$ are quite simple;
beginning again with \eqref{eq:inductivestep}, we obtain $n+1$
terms, all of which are rather simple for $n=2,3$. The details are
left to the reader.

Consider $n=5$, the most complicated case remaining.  From \eqref{eq:inductivestep} 
we obtain an upper bound of
\begin{equation}  \label{eq:casefivemainpart}
C_5 b_5 \tilde g^1_1 (g^1_1)^4
+ C_3 b_2\tilde g^1_1(g^1_1)^2(g^1_2)^2
+ C_2 b_3\tilde g^1_1 g^1_1(g^1_2)^3
+ b_5  C_5\tilde g^1_2(g^1_2)^4
\end{equation}
plus the two special terms involving $\int_{E^1_1} \tilde f$,
$\int_{E^1_2}f$. 
By Cauchy-Schwarz we may bound the square of \eqref{eq:casefivemainpart}
by 
\begin{equation*}
b_5(C_5\tilde g^1)^2
\big[(g^1_1)^4 + \beta(g^1_1)^2(g^1_2)^2 + \beta g^1_1(g^1_2)^3\big]^2
+ (g^1_2)^4,
\end{equation*}
where $\beta$ may be made as small as desired by choosing $C_5$
to be sufficiently large relative to $C_2,C_3$.
Thus the desired inequality reduces to
\begin{equation*}
(x^4+\beta x^2y^2+\beta xy^3)^2+y^8\le (x^2+y^2)^4,
\end{equation*}
which holds for small $\beta$.
The case $n=4$ is similar but simpler since one fewer term
appears in the analogue of \eqref{eq:casefivemainpart}.

We next discuss $M_n^*$. 
Let $\chi_y$ denote the characteristic function of $(-\infty,y]$
and apply the result just proved to $M_n(\tilde f,f\cdot\chi_y,
\cdots,f\cdot\chi_y)$ to obtain
\begin{equation}
M_n^*(\tilde f,f,\dots,f)
\le \frac{C^{n+1}}{\sqrt{n!}}
g_{-\delta}(\tilde f) \cdot[\sup_y g_\delta(f\cdot\chi_y)]^{n-1}.
\end{equation}
It was shown in the proof of Proposition~4.2 of 
\cite{christkiselevfiltrations} that 
for any function $F$,
\begin{equation}
\sup_y 
\sum_{m=1}^\infty
(\sum_j |{\textstyle\int}_{E^m_j} F\chi_y|^2)^{1/2}
\le C
\sum_{m=1}^\infty
m (\sum_j |{\textstyle\int}_{E^m_j} F|^2)^{1/2},
\end{equation}
and we may dominate the multiplicative coefficient
$m$ by $C_\epsilon 2^{m\epsilon}$ on the right-hand side.
Exactly the same proof allows the insertion of a factor of $2^{\delta m}$
after the first sum on both sides of the inequality, 
provided that $\delta\ge 0$, yielding the bound asserted.

We have so far
discussed only the special case where $f_2=f_3=\cdots = f_n =f$,
but the general case is treated by exactly the same argument, 
simply bounding
$|\int_{E^m_j}f_i|$ by $\max_k |\int_{E^m_j}f_k|$ wherever the former
arises in the proof.
\end{proof}

\section{Multilinear operators} \label{section:multilinear}

The bounds derived in the preceding section
for multiple integrals will be used
to obtain Lebesgue space norm bounds for certain multilinear operators.
Consider a family of integral operators
\begin{equation}
T_if(\lambda) = \int_\reals K_i(\lambda,x)f(x)\,dx.
\end{equation}
Denote by $\|T\|_{p,q}$ the norm of $T$ as an operator from
$L^p(\reals)$ to $L^q(\reals)$.
Consider associated multilinear operators
\begin{equation}
\scriptt_n(f_1,f_2,\dots,f_n)(\lambda)
= \iint_{x_1\le x_2\le\cdots\le x_n}
\prod_{i=1}^n K_i(\lambda,x_i)f_i(x_i)dx_i
\end{equation}
and their maximal variants
\begin{equation}
\scriptt_n^*(f_1,f_2,\dots,f_n)(\lambda)
= \sup_y\big|
\iint_{x_1\le x_2\le\cdots\le x_n\le y}
\prod_{i=1}^n K_i(\lambda,x_i)f_i(x_i)dx_i
\big|\ .
\end{equation}

\begin{theorem} \label{thm:multilinearoperators}
For any $N<\infty$, $p<q$,  and $2\le q$,
there exists $C<\infty$ such that for any $n\ge 2$ and any collection
of operators $\{T_i: 1\le i\le n\}$ of cardinality $\le N$,
for any collection of functions $\{f_1,\dots,f_n\}$ of cardinality
$\le N$,
\begin{equation}
\|\scriptt_n^*(f_1,\dots,f_n)\|_{L^r}
\le \frac{C^{n+1}}{\sqrt{n!}}
\|T_1\|_{2,2}\|f_1\|_{L^2}\prod_{i=2}^n\big(\|T_i\|_{p,q}
\|f_i\|_{L^p}
\big)
\end{equation}
where
\begin{equation}
\frac1r = \frac12 + \frac{n-1}q\ .
\end{equation}
\end{theorem}

\begin{proof}
It suffices to prove this under the assumption that
$\|f_1\|_{L^2}=1=\|f_i\|_{L^p}$ for every $i\ge 2$.
Construct a martingale structure $\{E^m_j\}$ for $\reals$ so that 
\begin{equation}
\int_{E^m_j}|f_1|^2
+ \sum_{i\ge 2}^* \int_{E^m_j}|f_i|^p
= 2^{-m} 
\int_\reals 
\big(
|f_1|^2
+ \sum_{i\ge 2}^*|f_i|^p\big)
\end{equation}
where the notation $\sum\limits^*$ indicates that the sum is to be taken 
over a maximal set of indices $i$ for which the functions $f_i$
are distinct.
Denote by $\chi^m_j$ the characteristic function of $E^m_j$.

Let $\delta>0$ be a constant to be specified. 
By Proposition~\ref{prop:numericalbound},
\begin{equation}  \label{eq:boundfromtheorem1}
\scriptt_n^*(f_1,\dots,f_n)(\lambda)
\le \frac{C^{n+1}}{\sqrt{n!}}
\tilde G(\lambda)
G(\lambda)^{n-1}
\end{equation}
where
\begin{align}
\tilde G(\lambda)
& = \sum_{m=1}^\infty 2^{-\delta m} 
(\sum_j |T_1(f_1\cdot \chi^m_j)(\lambda)|^2)^{1/2}
\\
G(\lambda)
& = \sum_{m=1}^\infty 2^{\delta m} 
(\sum_j \max_i|T_i(f_i\cdot \chi^m_j)(\lambda)|^2)^{1/2}.
\end{align}
Now 
\begin{multline}
\|\tilde G\|_{L^2} \le \sum_m 2^{-\delta m}
(\sum_j\|T_1(f_1\cdot\chi^m_j)\|_{L^2}^2)^{1/2}
\\
\le \sum_m 2^{-\delta m}
\|T_1\|_{2,2}(\sum_j\|f_1\cdot\chi^m_j)\|_{L^2}^2)^{1/2}
\le C_\delta
\|T_1\|_{2,2}\|f_1\|_{L^2}
\end{multline}
provided that $\delta>0$.
As for $G$, we may majorize 
\begin{equation}
G(\lambda)
\le \sum^*_i  \sum_{m=1}^\infty 2^{\delta m} 
(\sum_j |T_i(f_i\cdot \chi^m_j)(\lambda)|^2)^{1/2}.
\end{equation}
It is shown on page 413 of \cite{christkiselevfiltrations}
(see also the proof of Corollary~\ref{cor:nontangentialmaxbound} below)
that since $p<2\le q$, there exists $\e>0$ such that
\begin{equation}
\|(\sum_j |T_i(f_i\cdot \chi^m_j)|^2)^{1/2}\|_{L^q}
\le C2^{-\e m},
\end{equation}
under the condition that $\|f_i\chi^m_j\|_{L^p}^p\le C' 2^{-m}$
for all $j,m$.
Choosing $\delta<\e$ and summing over $m$ gives
\begin{equation}
\|G\|_{L^q} \le C<\infty,
\end{equation}
of course under the hypothesis that $\|f_i\|_{L^p}=1$ for all $i\ge 2$.
An application of H\"older's inequality concludes the proof.
\end{proof}
 
Our main theorems are based on estimates for such multilinear operators;
however, we require not the conclusion of 
Theorem~\ref{thm:multilinearoperators}, but rather the
more detailed information contained in \eqref{eq:boundfromtheorem1}
together with the norm bounds for $\tilde G,G$.

\section{Nontangential convergence and the maximal function}\label{nontan}

Throughout the paper we write $\complex^+=\{z\in\complex: \Im(z)>0\}$.
Consider the cones 
\begin{align*}
\Gamma_{\alpha}(w) &=\{z\in\complex^+: |\Re(z)-w|<\alpha\Im(z)\},
\\
\Gamma_{\alpha,\delta}(w)
&=\{z\in\complex^+\cap\Gamma_\alpha(w): \Im(z)<\delta\}.
\end{align*}
A function $f(z)$ defined on $\complex^+$
is said to converge to a limit $a$ as $z\to w$ nontangentially 
if for every $\alpha<\infty$,
$f(z)\to a$ as $z\to w$ within the cone $\Gamma_\alpha(w)$.
The nontangential maximal function of $f$ is defined by
\begin{equation}
Nf(w) =N_{\alpha,\delta}f(w)= \sup_{z\in \Gamma_{\alpha,\delta}(w)}|f(z)|.
\end{equation}

We will need the following local variant of a standard property of 
functions holomorphic in the whole half plane $\complex^+$.
\begin{lemma}\label{analyticf}
Let $\delta>0$. Let $\Lambda$ be an open subinterval of $\reals$,
and let $\Lambda'\Subset\Lambda$ be a relatively compact subinterval.
Let $0\le q\le\infty$. Let $B$ be a Banach space.
Suppose that $F$ is a holomorphic $B$-valued function
in  $\{\lambda+i\e: 0<\e<\delta,\ \lambda\in\Lambda\}$.
Then for any $\alpha<\infty$,
there exist $C<\infty$ and $\delta'>0$
depending on $\alpha,\Lambda,\Lambda'$, but not on $F$, such that
\begin{equation}
\|N_{\alpha,\delta'}F(\cdot)\|_{L^q(\Lambda')} \le
C \sup_{0<\e<\delta}
\|F(\cdot+i\e)\|_{L^q(\Lambda)} .
\end{equation}
\end{lemma}

By a $B$-valued holomorphic function we mean a continuous function $f$ from
an open subset of $\complex$ to $B$, such that $z\mapsto \ell(f)(z)$
is holomorphic for every bounded linear functional $\ell$ on $B$.

\begin{proof}
Suppose first that $f$ is continuous on the closed rectangle
$\Lambda+i[0,\delta]$, and that $1<q<\infty$. 
Let $0<\delta''<\delta$ be a small number to be chosen.
Choose a third interval $\Lambda''$, so that $\Lambda'\subset\Lambda''
\subset\Lambda$ and the distance from each interval to the boundary of the next
larger interval is strictly positive.
Consider the rectangle $R = \Lambda''+i[0,\delta'']=\{w: \Re(w)\in R
\text{ and } \Im(w)\in[0,\delta'']\}$.  

For any $z\in R$ with $\Re(z)\in\Lambda'$ and $0<\Im(z)<\delta''/2$,
we may write $f(z) = \int_{\p R}f(\zeta)d\omega_z(\zeta)$ 
where $\omega_z$ is harmonic measure
on $\p R$.
Fix a conformal map $\varphi$ from $R$ to the unit disk; this map is smooth everywhere except
at the corners of $R$, where it behaves locally like $z\mapsto z^2$ in appropriate
local coordinates.
Harmonic measure for $R$ may be computed by pulling back
harmonic measure from the unit disk under $\varphi$.  From 
this we readily deduce that for any $z\in R$ with $\Re(z)\in\Lambda$,
\begin{equation} \label{poissonrepresentation}
f(z) = \int_{\Lambda''} k(z,\zeta)f(\zeta)\,d\sigma(\zeta) 
+ O\Big(\int_{\p R\backslash\Lambda''} |f(\zeta)|\,d\sigma(\zeta)\Big)
\end{equation}
where $\sigma$ denotes arc length measure on $\p R$, and where
$k$ satisfies upper bounds of Poisson kernel type:
$|k(x+it,y)|\le Ct[|x-y|^2+t^2]^{-1}$.

Denote by $F$ the restriction of $f$ to $\Lambda''$,
and by $M$ the Hardy-Littlewood maximal function.
We combine these bounds with
the usual majorization of the nontangential maximal
function of a Poisson integral by $M$
to conclude that for $z\in \Gamma_{\alpha,\delta'}(y)$ with $0<\Im(z)<\delta''$
and $y\in\Lambda'$, the first term on the right-hand side of 
\eqref{poissonrepresentation} is bounded by $C'MF(y)$. 

The second term on the right is not suitably bounded, in general. Therefore we
consider all translates $R_s = \{z+s: z\in R\}$ where $s\in \reals$ ranges over
a small interval $I$ centered at $0$. For each $R_s$ we have a variant of
\eqref{poissonrepresentation}, obtained by conjugating with the translation
$z\mapsto z+s$; the integral in the first term now extends over $\Lambda''+s$.
By averaging all these variants over $s\in I$ with respect to normalized Lebesgue
measure, we obtain
\begin{equation*}
f(z) = \int_{\tilde\Lambda} \tilde k(z,w)d\sigma(w) + O\big(\int_{Q} |f(w)|\,dw\big)
+ O\big( \int_{\tilde\Lambda} |f(w+i\delta'')|\,d\sigma(w) \big)
\end{equation*}
where $\tilde\Lambda=\cup_{s\in I} \Lambda''+s $ is slightly larger than $\Lambda''$,
$\tilde k$ satisfies the same bounds as $k$, $dw$ denotes Lebesgue measure in
$\complex$, and $Q$ is a certain compact subset of the closed rectangle 
$\Lambda+i[0,\delta]$.
The nontangential maximal function of the
first term is majorized by $CMF$, just as before. The second is majorized, uniformly
in $z$, by $\sup_\e \|f(\cdot+i\e)\|_{L^q(\Lambda)}$.
The third is already under control, by hypothesis, since $q\ge 1$.
Because $M$ is bounded on $L^q$,
we obtain the desired conclusion, under the supplementary hypothesis that $f$
extends continuously to $\Im(z)=0$.

For general $f$ not necessarily continuous up to $\Im(z)=0$, 
consider $f_\eta(z) = f(z+i\eta)$, for small $\eta\in\reals^+$.
Apply the result just proved to $f_\eta$, and pass to the limit $\eta\to 0^+$
using Fatou's lemma.

In order to extend this argument to $0<q\le 1$, it suffices to recall that 
if $F$ is a $\complex$-valued holomorphic function, then $|F|^s$ is subharmonic
for any $s>0$. Fixing any $0<s<q$ any bounded linear functional $\ell$ on $B$,
and setting $F(z) = \ell(f(z))$, 
we may therefore write instead of \eqref{poissonrepresentation} 
\begin{equation*} 
|F(z)|^s \le \int_{\Lambda''} k(z,\zeta)|F(\zeta)|^s\,d\sigma(\zeta)
+ O\Big(\int_{\p R\backslash\Lambda''} |F(\zeta)|^s\,d\sigma(\zeta)\Big).
\end{equation*}
The proof then proceeds as above, using the fact that $f\mapsto [M(|f|^s)]^{1/s}$
is bounded on $L^q$ for all $q>s$.
\end{proof}

\section{Generalized eigenfunctions associated to complex 
spectral parameters}  \label{section:complex}


Consider the generalized eigenfunction equation
\begin{equation} \label{eigenfneqn}
-u'' + V(x)u = zu,
\end{equation}
where $z$ is permitted to be complex.
In earlier work we have analyzed the solutions to this equation for $z$ 
real and positive, and have shown that for almost every such $z$ there 
exists a solution with certain asymptotic behavior as $x\to+\infty$. 
Our present purpose is to analyze solutions for complex $z$ and to show 
that they tend almost everywhere to the solutions for real $z$, 
as $z$ approaches the positive real axis.  

We restrict attention to parameters $z\in\complex^+\cup \reals^+$,
and choose a branch of $\sqrt{z}$ which has nonnegative imaginary part
for such $z$.
Define the phases
\begin{equation}
\xi(x,z) = \sqrt{z}x-(2\sqrt{z})\rp\int_0^x V.
\end{equation}

\begin{theorem}  \label{thm:complexeigenfunctions}
Let $1\le p<2$ and assume that $V\in L^1+L^p$.
For each $z\in\complex^+$ there exists a (unique) solution $u(x,z)$ of  
the generalized eigenfunction equation \eqref{eigenfneqn}
satisfying
\begin{equation} \label{wkbasymptotics}
u(x,z) - e^{i\xi(x,z)} \to 0 
\text{ and }
\partial u(x,z)/\partial x - i\sqrt{z}e^{i\xi(x,z)} \to 0
\qquad\text{as } x\to +\infty.
\end{equation}
$u(x,z)$ is continuous as a function on $\complex^+\times\reals$,
and is holomorphic with respect to $z$ for each fixed $x$.

Likewise, there exists such a (unique) solution
for almost every $z\in\reals^+$.
For almost every $E\in\reals$,
$u(x,z)$ converges to $u(x,E)$ uniformly for all $x$ in any 
interval bounded below, as $z\to E$ nontangentially.
\end{theorem}

Here ``almost every'' means with respect to Lebesgue measure.
The existence of such solutions for almost every $z\in\reals^+$ is proved
in \cite{christkiselevdecaying}. For $z\in\complex^+$ it is well known
under weaker hypotheses on $V$. 
The new point here is the convergence as $z\to\reals^+$, and in particular, 
the fact that it is globally uniform in $x$.

It suffices to prove the theorem for $x\in [0,\infty)$, since the
conclusion for $x\in [\rho,\infty)$, for any $\rho>-\infty$,
then follows via the eigenfunction equation \eqref{eigenfneqn}.

By rewriting \eqref{eigenfneqn} as a first-order system, performing 
a couple of algebraic transformations, reducing to an integral equation,
and solving it by iteration,
one arrives \cite{christkiselevdecaying} at a formal 
series representation for solutions of \eqref{eigenfneqn}:
\begin{equation}  \label{formalsolutionseries}
\begin{pmatrix} u(x,z) \\ u'(x,z) \end{pmatrix}
=
\begin{pmatrix}
e^{i\xi(x,z)} & e^{-i\xi(x,z)} \\
i\sqrt{z} e^{i\xi(x,z)}   & -i\sqrt{z}  e^{-i\xi(x,z)}
\end{pmatrix}
\cdot
\begin{pmatrix}
\sum_{n=0}^\infty T_{2n} (V, \dots, V)(x,z)
 \\
-\sum_{n=0}^\infty T_{2n+1} (V, \dots, V)(x,z) 
\end{pmatrix}
\end{equation}
where
\begin{equation}
 T_n (V_1, \dots, V_n)(x,z)= (2\sqrt{z})^{-n}
\int_{x \leq t_1 \leq \dots \leq t_n} \prod_{j=1}^n
e^{2i (-1)^{n-j} \xi(t_j,z)} V_j(t_j)\,dt_j,
\end{equation}
with the convention $T_0(\cdot) \equiv 1$.
To prove Theorem~\ref{thm:complexeigenfunctions}
we will show that each multilinear expression $T_n$ 
is well-defined for all $z\in\complex^+$, 
that $T_n(\cdot)(x,z)\to 0$ as $x\to+\infty$ for all $n\ge 1$,
that they have the natural limits as $z\to\reals^+$ nontangentially,
and that these expressions satisfy bounds sufficiently strong to enable
us to sum the infinite series to obtain the desired conclusions.

Substitute $\zeta=\sqrt{z}$ and write $\zeta = \lambda+i\e$,
noting that $\e>0$.
Also write
$\phi(x,\zeta) = \xi(x,z)$, $S_n(V_1,\dots)(x,\zeta)
= T_n(V_1,\dots)(x,z)$.
Let $\{E^m_j\}$ be a martingale structure on $\reals^+$.
Denote by $t_{m,j}^\pm$ respectively the right ($+$) and left ($-$)
endpoints of the interval $E^m_j$. 

The real part of the exponent 
$2i\sum_{j=1}^n (-1)^{n-j} \xi(t_j,z)$
is to leading order
\begin{equation*}
-2\Im(\sqrt{z}[(t_n-t_{n-1})+(t_{n-2}-t_{n-3})+\cdots])
= -2\e\cdot [(t_n-t_{n-1})+(t_{n-2}-t_{n-3})+\cdots],
\end{equation*}
which is nonnegative (for $x\ge 0$) 
for all $z\in\complex^+$ since $t_1\le t_2\cdots$; the exponential
factor decays rapidly as 
$[(t_n-t_{n-1})+(t_{n-2}-t_{n-3})+\cdots]\to\infty$.
This accounts for the difference between $z\in\complex^+$
and $z\in\reals^+$.

Define
\begin{equation}\begin{split}
G_m(V)(\zeta)
&= 
\Big(
\sum_{j=1}^{2^m}
|s^{m,-}_j(V,\zeta)|^2 + | s^{m,+}_j(V,\zeta)|^2\Big)^{1/2}
\\
G(V) &= \sum_{m=1}^\infty mG_m(V)
\end{split}\end{equation}
where
\begin{equation}  \label{emjintegrals} \begin{split}
s^{m,-}_j(V,\zeta)
&= \int_{E^m_j} e^{2i[\phi(t,\zeta)-\phi(t^-_{m,j},\zeta)]}V(t)\,dt
\\
s^{m,+}_j(V,\zeta) &= \int_{E^m_j} e^{2i[\phi(t^+_{m,j},\zeta)-\phi(t,\zeta)]}V(t)\,dt.
\end{split}\end{equation}
When $j=2^m$, and only then, the right endpoint of $E^m_j$ is infinite. 
To simplify notation we make the convention that
for $j=2^m$, the second term
$|s^{m,+}_j(V,\zeta)|^2$
is always to be omitted, in the definition of $G$ and anywhere
else that the quantities $\phi(t^+_{m,j},\zeta)$ arise.

The following definitions will allow us to regard $G$ as a linear operator,
and hence to exploit properties of holomorphic functions.
\begin{definition}
$\scriptb$ denotes the Banach space consisting of all 
sequences $\{\complex^2\owns s^m_j: m\ge 0,\ 1\le j\le 2^m\}$, with the norm
$\|s\|_\scriptb = \sum_m m\big(\sum_{j=1}^{2^m} |s^m_j|^2\big)^{1/2}$.

$\frakG: L^p\mapsto\scriptb$ denotes the operator
\begin{equation}
\frakG(V)(\zeta) = \{(s^{m,+}_j(V,\zeta),s^{m,-}_j)(V,\zeta): 1\le m<\infty, 1\le j\le 2^m\}.
\end{equation}
\end{definition}
Thus
\begin{equation*}
G(V)(\zeta) = \|\frakG(V)(\zeta)\|_{\scriptb}.
\end{equation*}
Likewise we may write $G_M(V) = \|\frakG_M(V)\|_{\scriptb}$
with the analogous definition of $\frakG_M$.

The real parts of $i\phi(t,\zeta)-i\phi(t^-_{m,j},\zeta)$
and $i\phi(t^+_{m,j},\zeta)-i\phi(t,\zeta)$
are bounded above uniformly for $0<\Im(\zeta)\le 1$,
and are $\le -c\Im(\zeta)(t-t^-_{m,j})$ and
$\le -c\Im(\zeta)(t^+_{m,j}-t)$, respectively;
see \eqref{realpartnonpositive}.
Therefore for $V\in L^1+L^\infty$ and $\zeta\in\complex^+$,
each of these integrals converges absolutely, and each defines a holomorphic
scalar-valued function of $\zeta\in\complex^+$.
Thus $\frakG(V)$ may be regarded as a $\scriptb$--valued holomorphic function\footnote{
By this we mean simply that it is a continuous mapping into $\scriptb$
with respect to the norm topology, and that each $s^{m,\pm}_j$ is
a holomorphic scalar-valued function.}
in any open set where it can be established that the series 
defining $\|\frakG(V)\|_\scriptb$ converges uniformly.

We will also use the following
variant. Given a collection of functions $(V_1,\dots,V_n)$, we define
\begin{equation}
G(\{V_i\})(\zeta)
= \sum_{m=1}^\infty m\cdot
\Big(
\sum_{j=1}^{2^m}
\sum^*_i
|s^{m,+}_j(V_i,\zeta)|^2+|s^{m,-}_j(V_i,\zeta)|^2
\Big)^{1/2},
\end{equation}
where $\sum\limits^*$ indicates that the sum is taken over a maximal 
set of indices $i$ for which the functions $V_i$ are all distinct.

\begin{lemma} \label{lemma:Gmajorizationwithmodifiedphases}
For all $V,n$ and all $\zeta\in \complex^+\cup\reals^+$,
\begin{equation}
\sup_{x\in\reals^+}|S_n(V,V,\dots,V)(x,\zeta)|
\le \frac{C^{n+1}}{\sqrt{n!}} 
G(V)(\zeta)^n.
\end{equation}
More generally,
for all $n$ and all $\{V_1,\dots,V_n\}$,
\begin{equation}
\sup_{x\in\reals^+}|S_n(V_1,V_2,\dots,V_n)(x,\zeta)|
\le \frac{C_k^{n+1}}{\sqrt{n!}} 
G(\{V_i\})(\zeta)^n,
\end{equation}
provided that $\{V_i\}_{i=1}^n$ has cardinality $\le k$. 
\end{lemma}

Write $p'=p/(p-1)$.
\begin{lemma} \label{lemma:parseval}
For any compact interval $\Lambda\Subset(0,\infty)$,
there exists $C<\infty$ such that
for any $1\le p\le 2$, 
for all $t'\in\reals$ and $f\in L^p(\reals)$,
for every $\e\ge 0$,
\begin{align}
&\|\int_{t\ge t'} e^{2i[\phi(t, \lambda+i\e) - \phi(t', \lambda+i\e)]}
f(t)\,dt\|_{L^{p'}(\Lambda,d\lambda)}
\le C\|f\|_{L^p}
\\
&\|\int_{t\le t'} e^{2i[\phi(t', \lambda+i\e) - \phi(t, \lambda+i\e)]}
f(t)\,dt\|_{L^{p'}(\Lambda,d\lambda)}
\le C\|f\|_{L^p}.
\end{align}
\end{lemma}

For $\e=0$ these integrals need not converge absolutely, hence require
interpretation. They are initially well-defined for compactly 
supported $f$, and the lemma asserts an {\em a priori} bound for 
such functions.  Then they are defined for general $f\in L^p$ 
by approximating in $L^p$ norm by compactly supported functions, 
and passing to the limit in $L^{p'}$ norm.

Let an exponent $p<\infty$ be specified.
Recall that a martingale structure is said to be adapted to $f$ in $L^p$
if $\int_{E^m_j}|f|^p = 2^{-m}\int|f|^p$ for all $m,j$.
Recall from the preceding section the definitions of the cones $\Gamma_{\alpha,\delta}$
and associated nontangential maximal functions $N_{\alpha,\delta}$.

\begin{corollary} \label{cor:nontangentialmaxbound}
Let $\alpha<\infty$, let $1\le p \le 2$, and let
$\Lambda\Subset(0,\infty)$ be any compact subinterval.
Then there exist $C<\infty,\delta>0$ such that for any
$f\in L^p(\reals)$ and for any martingale structure
$\{E^m_j\}$ on $\reals^+$ 
\begin{equation}  \label{nonadaptedbound}
\|N_{\alpha,\delta} G_m(f)(\lambda)\|_{L^{p'}(\Lambda,d\lambda)}
\le C\|f\|_{L^p}.
\end{equation}
Moreover for each $1\le p<2$ 
there exists $\rho>0$ such that for any $f\in L^p$ and
for any martingale structure adapted
to $f$ in $L^p$, 
\begin{equation} \label{GMdecaybound}
\|N_{\alpha,\delta} G_m(f)(\lambda)\|_{L^{p'}(\Lambda,d\lambda)}
\le C2^{-\rho m}\|f\|_{L^p}.
\end{equation}

Consequently under these additional hypotheses,
\begin{equation}
\|N_{\alpha,\delta} G(f)(\lambda)\|_{L^{p'}(\Lambda,d\lambda)}
\le C\|f\|_{L^p}.
\end{equation}
Moreover, for almost every $\lambda\in\Lambda$,
\begin{equation}
\|\frakG(f)(\zeta)-\frakG(f)(\lambda)\|_{\scriptb} \to 0
\end{equation}
as $\zeta\to\lambda$ nontangentially.
\end{corollary}

We postpone the proofs of Lemmas~\ref{lemma:Gmajorizationwithmodifiedphases}, 
\ref{lemma:parseval} and Corollary~\ref{cor:nontangentialmaxbound}
until the end of the section.

For any $t\ge t'$ and any $\zeta=\lambda+i\e$ with $\e\ge 0$,
\begin{equation} \label{realpartnonpositive}
\Re(i\phi(t,\lambda+i\e)-i\phi(t',\lambda+i\e))
= -\e(t-t') -\e(\lambda^2+\e^2)\int_{t'}^t V,
\end{equation}
and $|\int_{t'}^t V|\le |t-t'|^{1/p'}\|V\|_{L^p}$.
Therefore
\begin{equation}
|e^{i[\phi(t,\lambda+i\e)-\phi(t',\lambda+i\e)]}|
\le Ce^{-c\e|t-t'|}
\end{equation}
where $C,c\in\reals^+$ are constants which depend only on the $L^p$ norm
of $V$.
Hence for all $n\ge 1$,
\begin{align}
|S_{2n}(V,V,\dots,V)(x,\lambda+i\e)|
& \le \iint_{x\le t_1\le\cdots\le t_{2n}} 
e^{-c\e(t_{2n}-t_{2n-1}+t_{2n-2}-\cdots)}
\prod_j V(t_j)\,dt_j
\notag
\\
& \le \frac{C^n}{n!} 
\Big(\iint_{x\le t'\le t} e^{-c\e(t-t')}|V(t')V(t)|\,dt'\,dt\Big)^n
\notag
\\
& \le \frac{C^n}{n!} 
\e^{-2n(1-p\rp)} \|V\|_{L^p([x,\infty)}^{2n}.
\end{align}
In the same way, for $x\ge 0$, one obtains for $n\ge 0$
\begin{equation}
|S_{2n+1}(V,V,\dots,V)(x,\lambda+i\e)|
\le \frac{C^n}{n!} 
\e^{-(2n+1)(1-p\rp)}  \|V\|_{L^p([x,\infty)}^{2n+1}.
\end{equation}

Let $z$ belong to any compact subset of $\complex^+$.
Then $\e=\Im(\sqrt{z})$ has a strictly positive lower bound.
Therefore the individual terms of the series \eqref{formalsolutionseries} 
defining a formal solution of \eqref{eigenfneqn} do define uniformly bounded
functions of $(x,\zeta)$ for $x\ge 0$ and $\zeta$ in any compact subset
of $\complex^+$, and moreover, the 
series are uniformly absolutely convergent.
As in Lemma~4.2 of \cite{christkiselevdecaying}, 
it follows that the sums of these series define solutions of the 
generalized eigenfunction equation \eqref{eigenfneqn} for all such $z$.
Because the $L^p$ norm of $V$ over $[x,\infty)$ tends to zero 
as $x\to+\infty$, only $T_0(x)$ contributes in that limit,
so these solutions do have the desired WKB asymptotics $\exp(i\xi(x,z))$.
Clearly each summand depends holomorphically on $\zeta$, hence so do
the sums.

Existence of a (unique) solution for almost every $z\in\reals^+$
is proved in \cite{christkiselevdecaying}. 
We come now to the main step, where we relate 
$z\in\complex^+$ to $z\in\reals^+$.
For compactly supported $f\in L^1$,
the quantities $s^{m,\pm}_j(f,\zeta)$
are clearly holomorphic functions of $\zeta$ where
$\Im(\zeta)>0$, and are continuous at $\e=0$ for each $0\ne \lambda\in\reals$.
The same holds for $S_n(f,f,\dots,f)(x,\zeta)$, for each $x\in\reals$.

$\frakG(f)$ is likewise a $\scriptb$--valued holomorphic function 
of $\zeta$, continuous on $\complex^+\cup\reals^+$,
for compactly supported $f\in L^1$.
This follows by combining the holomorphy of the individual terms
\eqref{emjintegrals} with the rapid convergence bound \eqref{GMdecaybound}.

\begin{lemma}  \label{lemma:Snconvergence}
Let $1\le p<2$ and $V\in L^p(\reals)$.
For almost every $E\in\reals$,
for every $n\ge 1$,
$T_n(V,V,\dots,V)(x,z)\to T_n(V,V,\dots,V)(x,E)$
uniformly for all $x\ge 0$
as $\complex^+\owns z\to E$ nontangentially.
\end{lemma}

\begin{proof}
It is equivalent to show that 
\begin{equation*}
\sup_{x\ge 0}\big|
S_n(V,V,\dots,V)(x,\zeta)-S_n(V,V,\dots,V)(x,\lambda)\big| \to 0
\end{equation*}
as $\zeta\to\lambda$ nontangentially, for almost all $\lambda\in\Lambda$,
for any fixed compact interval $\Lambda\subset\reals^+$.
For $V=W\in L^1$ with compact support, we have already established
convergence uniformly in $x,\lambda$, as $\complex^+\owns\zeta\to \lambda$
unrestrictedly (rather than merely nontangentially),
since the phases $\phi(t,\zeta)$ converge uniformly to $\phi(t,\lambda)$
for $t,\zeta,\lambda$ in any compact set.

Let $V$ be given, and remain fixed for the remainder of this proof.
Set $\phi_V(t,\zeta)=\zeta t - (2\zeta)\rp\int_0^t V$.
Whenever we write $S_n(f_1,f_2,\dots,f_n)$, it is defined
in terms of the phases $\phi(t_i,\zeta) = \phi_V(t_i,\zeta)$,
independent of $f_1,f_2,\dots$. Thus the $S_n$ are here
genuine multilinear operators.

Let $\e>0$ be arbitrary, and fix a martingale structure $\{E^m_j\}$
adapted to $V$ in $L^p$ on $\reals^+$.
Decompose $V=W+(V-W)$
where $W(x) =V(x)\chi_{(0,R]}(x)$, with $R$ chosen so
that $\|V-W\|_{L^p}<\e$
and moreover so that $\|NG(V-W)\|_{L^{p'}(\Lambda)}<\e$.
Such a choice is possible, since 
\begin{equation*}
\|NG_M(V\chi_{[R,\infty)})\|_{L^{p'}(\Lambda)}
\le C\min(2^{-M\delta}\|V\|_{L^p},\,\| V\chi_R\|_{L^p}).
\end{equation*}
Then
\begin{multline}
|S_n(V,V,\dots,V)(x,\zeta)-S_n(V,V,\dots,V)(x,\lambda)|
\\
\le
|S_n(W,W,\dots,W)(x,\zeta)-S_n(W,W,\dots,W)(x,\lambda)|
\\
+
|S_n(V,V,\dots,V)(x,\zeta)-S_n(W,W,\dots,W)(x,\zeta)|
\\
+
|S_n(V,V,\dots,V)(x,\lambda)-S_n(W,W,\dots,W)(x,\lambda)| .
\end{multline}
The first term on the right tends to zero, in the sense desired.
Majorize the second by
\begin{multline}
|S_n(V,V,\dots,V)(x,\zeta)-S_n(W,W,\dots,W)(x,\zeta)|
\\
\le
\sum_{i=1}^n
|S_n(V,\dots,V,V-W,W,\dots,W)(z,\zeta)|
\end{multline}
where in the $i$-th summand, the argument of $S_n$ has
$i-1$ copies of $V$ and $n-i$ copies of $W$.
Fix any aperture $\alpha\in\reals^+$.
Thus as established in the proof of
Proposition~4.1 of \cite{christkiselevdecaying},
\begin{multline} 
\sup_{x\ge 0}
\sup_{\zeta\in \Gamma_{\alpha,\delta} (\lambda)}
|S_n(V,V,\dots,V)(x,\zeta)-S_n(W,W,\dots,W)(x,\zeta)|
\\
\le
C_n
\sum_{i=1}^n
\sup_{\zeta\in \Gamma_{\alpha,\delta}(\lambda)}
G(V)^{i-1}(\zeta)G(W)^{n-i}(\zeta)G(V-W)(\zeta)
\\
\le 
C_n
\sum_{i=1}^n
NG(V)^{i-1}(\lambda)NG(W)^{n-i}(\lambda)NG(V-W)(\lambda).
\end{multline}

Let $q=p'$.
By Chebyshev's inequality,
for any $\beta>0$,
\begin{equation}\begin{split}
&|\{\lambda\in\Lambda:
\sup_{\zeta\in \Gamma_{\alpha,\delta}(\lambda)}\sup_{x\ge 0}
|S_n(V,V,\dots,V)(x,\zeta)-S_n(W,W,\dots,W)(x,\zeta)|>\beta\}|
\\
&\qquad\qquad \le
C_n\beta^{-q/n} \sum_{i=1}^n
\| NG(V)^{i-1}NG(W)^{n-i}NG(V-W)\|_{L^{q/n}(\Lambda)}^{q/n}
\\
&\qquad\qquad \le
C_n\beta^{-q/n} \sum_{i=1}^n
(\|NG(V)\|_{L^q(\Lambda)}^{i-1}
\|NG(W)\|_{L^q(\Lambda)}^{n-i}
\|NG(V-W)\|_{L^q(\Lambda)})^{q/n}
\\
&\qquad\qquad \le
C_n\beta^{-q/n} \sum_{i=1}^n
(\|V\|_{L^p}^{n-1}\|V-W\|_{L^p})^{q/n}
\\
&\qquad\qquad \le
C_n\beta^{-q/n}
\|V\|_{L^p}^{q(n-1)/n}\e^{q/n}.
\end{split}\end{equation}
The same reasoning gives 
\begin{multline}
|\{\lambda\in\Lambda:
\sup_{x\ge 0}
\big|S_n(V,V,\dots,V)(x,\lambda)-S_n(W,W,\dots,W)(x,\lambda)|>\beta\}\big|
\\
\le
C_n\beta^{-q/n}
\|V\|_{L^p}^{q(n-1)/n}\e^{q/n}.
\end{multline}
Consequently
\begin{multline}
|\{\lambda\in\Lambda:
\limsup_{\Gamma_\alpha(\lambda)\owns\zeta\to\lambda}\,\,
\sup_{x\ge 0}\,
\big|S_n(V,V,\dots,V)(x,\zeta)-S_n(V,V,\dots,V)(x,\lambda)\big|
>\beta\}
\\
\le
C_n\beta^{-q/n}
\|V\|_{L^p}^{q(n-1)/n}\e^{q/n},
\end{multline}
for all $\beta,\e\in\reals^+$.
Letting $\e\to 0$, we conclude that the $\limsup$ vanishes for almost
every $\lambda$.
\end{proof}

It is now straightforward to sum the series to obtain the same
conclusion regarding convergence of $u(x,z)=u(x,\zeta^2)$
to $u(x,E)=u(x,\lambda^2)$:
\begin{equation}
|u(x,\zeta^2)-u(x,\lambda^2)|
\le
\sum_{n=0}^\infty 
|S_n(V,V,\dots,V)(x,\zeta)-S_n(V,V,\dots,V)(x,\lambda)|
\end{equation}
and
\begin{equation}\begin{split}
&\sup_{\zeta\in\Gamma_{\alpha,\delta}(\lambda)}\sup_{x\ge 0}
\sum_{n=M}^\infty
|S_n(V,V,\dots,V)(x,\zeta)-S_n(V,V,\dots,V)(x,\lambda)|
\\
&\qquad\qquad\le
\sup_{\zeta\in\Gamma_{\alpha,\delta}(\lambda)}\sup_{x\ge 0}
\sum_{n=M}^\infty
(|S_n(V,V,\dots,V)(x,\zeta)|+|S_n(V,V,\dots,V)(x,\lambda)|)
\\
&\qquad\qquad\le
\sum_{n=M}^\infty \frac{C^{n+1}}{\sqrt{n!}}
\big(NG(V)(\lambda)+G(V)(\lambda)\big)^n
\\
&\qquad\qquad\le
\frac{C^{M+1}}{\sqrt{M!}}\big(NG(V)(\lambda)+G(V)(\lambda)\big)^M
\sum_{k=0}^\infty
\frac{C^{k+1}}{\sqrt{k!}}
\big(NG(V)(\lambda)+G(V)(\lambda)\big)^k
\\
&\qquad\qquad\le
\frac{C^{M+1}}{\sqrt{M!}}\big(NG(V)(\lambda)+G(V)(\lambda)\big)^M
\exp(C \big(NG(V)(\lambda)+G(V)(\lambda)\big)^2).
\end{split}\end{equation}
For almost every $\lambda\in\reals$, 
$\big(NG(V)(\lambda)+G(V)(\lambda)\big)<\infty$,
and hence this expression tends to zero as $M\to\infty$.
Coupled with the convergence established for the individual terms
$S_n$ in Lemma~\ref{lemma:Snconvergence}, this completes the proof of
Theorem~\ref{thm:complexeigenfunctions}, modulo the proofs of 
Lemmas~\ref{lemma:parseval} and 
\ref{lemma:Gmajorizationwithmodifiedphases}, 
and Corollary~\ref{cor:nontangentialmaxbound}.
\hfill\qed

\begin{proof}[Proof of Lemma~\ref{lemma:parseval}]
The proofs of the two inequalities are essentially the same,
so we discuss only the first.
The exponent here is $i\Phi(t,\lambda+i\e) 
= i\phi(t,\lambda+i\e)- i\phi(t',\lambda+i\e)
= (i\lambda-\e)(t-t') - (i\lambda+\e)(\lambda^2+\e^2)\rp(\int_{t'}^t V)$.
Since $t\ge t'$ and $|\int_{t'}^t V| \le C+C|t-t'|^{1/2}$,
the real part of $\Phi$ is bounded above, uniformly for all
$\lambda$ in any compact subinterval $\Lambda\Subset\reals\backslash\{0\}$
and $\e\ge 0$. Thus $L^1$ is mapped boundedly to $L^\infty(\Lambda)$,
uniformly in $\e$; by interpolation it suffices to prove the $L^2$
estimate.

Fix any cutoff function $\eta\in C^\infty(\reals\backslash\{0\})$
and consider
\begin{equation}
\int \Big|\int_{t>t'} 
e^{i\Phi(t,\lambda+i\e)}f(t)\,dt\Big|^2\eta(\lambda)\,d\lambda
=
\iint_{s,t\ge t'} f(t)\bar f(s) K(t,s)\,dt\,ds
\end{equation}
where 
\begin{equation}
K(t,s) 
=\int e^{\Psi(t,s,\lambda+i\e)}
\eta(\lambda)
\,d\lambda 
\end{equation}
with
\begin{multline}
\Psi(t,s,\lambda+i\e) = {2i\Phi(t,\lambda+i\e)-2i\bar\Phi(s,\lambda+i\e)}
\\
= 2i\big[\lambda(t-s)-\lambda(\lambda^2+\e^2)\rp{\textstyle\int}_s^t V\big]
\\
-2\e\big[(t-t')+(s-t') 
+ (\lambda^2+\e^2)\rp({\textstyle\int}_{t'}^t V)
+ (\lambda^2+\e^2)\rp({\textstyle\int}_{t'}^s V)\big].
\end{multline}
We claim that $|K(t,s)|\le C(1+|s-t|)^{-2}$,
uniformly in $\e\ge 0$; this would suffice to imply the $L^2$ bound.
The integrand itself is bounded, uniformly in all parameters,
so it suffices to restrict attention to the case where $|s-t|\ge C_0$,
where $C_0$ is a sufficiently large constant.
In that case we integrate by parts,
integrating $\exp(2i[\lambda(t-s)-\lambda(\lambda^2+\e^2)\rp\int_s^t V])$,
and differentiating 
$\eta(\lambda)\cdot
\exp(-2\e[(t-t')+(s-t')
+ (\lambda^2+\e^2)\rp(\int_{t'}^t V)
+ (\lambda^2+\e^2)\rp(\int_{t'}^s V))$,
noting that $\partial [\lambda(t-s)-\lambda(\lambda^2+\e^2)\rp\int_s^t V]
/\partial\lambda\ge |s-t|/2$ provided $C_0$ is chosen to be
sufficiently large.
Thus we gain a factor of $(s-t)\rp$.
On the other hand, differentiating the other exponential with respect
to $\lambda$ brings in an unfavorable term 
$O(\e \int_{t'}^{\max(s,t)}|V|)$.
After two integrations by parts, the integrand is
\begin{equation}
O(|s-t|^{-2})\cdot e^{-2\e(s-t')-2\e(t-t')}
\cdot O( 1+\e^2( \int_{t'}^{\max(s,t)}|V| )^2 )
= 
O(|s-t|^{-2}),
\end{equation}
uniformly in $\e$.
\end{proof}

\begin{proof}[Proof of Lemma~\ref{lemma:Gmajorizationwithmodifiedphases}]
The new feature here is the introduction of the modifying factors
$\exp( \pm 2i\phi(t_{m,j}^\pm,\zeta) )$; without these, this is proved
in \cite{christkiselevfiltrations}. We will merely indicate
the modification needed in the argument, referring to
\cite{christkiselevfiltrations} for the rest.
Consider
\begin{equation}
\iint_{x\le t_1\le\cdots\le t_n}e^{2i[\phi(t_n,\zeta)-\phi(t_{n-1},\zeta)
+ \phi(t_{n-2},\zeta)-\cdots ]}
f(t_1)f(t_2)\cdots f(t_n)\,dt_1\cdots\,dt_n.
\end{equation}
Decompose the region of integration $\{t=(t_1,\cdots,t_n):
x\le t_1\le\cdots\le t_n\}$ as $\cup_{k=0}^n\Omega_k$
where $\Omega_k=\{t: x\le t_1\le\cdots\le t_k\le t^+_{1,1}
= t^-_{1,2}\le t_{k+1}\le\cdots\le t_n\}$.
The total integral is
\begin{equation}\begin{split}
&\sum_k\iint_{\Omega_k}
e^{2i[\phi(t_n,\zeta)-\phi(t_{n-1},\zeta)+\cdots\pm\phi(t_{k+1},\zeta)
\mp\phi(t^-_{1,2},\zeta)]}
\cdot
e^{2i[\pm\phi(t^+_{1,1}
\mp \phi(t_k,\zeta)\pm\phi(t_{k-1},\zeta)\mp\cdots]}
\prod_{j=1}^n f(t_j)\,dt_j
\\
& \qquad=
\sum_k
\Big(
\iint_{t^-_{1,2}\le t_{k+1}\le\cdots\le t_n}
e^{2i[\phi(t_n,\zeta)-\phi(t_{n-1},\zeta)+\cdots 
+(-1)^{n-k-1}\phi(t_{k+1},\zeta)]}
\prod_{j=k+1}^n f(t_j)\,dt_j
\Big)
\\
&\qquad\qquad\cdot
\Big(
\iint_{x\le t_1\le\cdots\le t_k\le t^+_{1,1} }
e^{2i(-1)^{n-k}[\phi(t_k,\zeta) - \phi(t_{k-1},\zeta) + \cdots
+ (-1)^{n-1}\phi(t_1,\zeta)]}
\prod_{j=1}^k f(t_j)\,dt_j
\Big).
\end{split}\end{equation}
For each $k$, each of the two factors on the right-hand side has
the same form as the multiple integral with which we began,
except that when $n-k$ is odd, 
an extra factor of $-1$ appears in the exponent in 
the integral with respect to $dt_k\cdots dt_1$; this minus sign
destroys the bounds we seek, as is clear from 
\eqref{realpartnonpositive}.
Therefore when $n-k$ is odd, we rewrite the corresponding term
as the modified product
\begin{multline*}
\Big(
\iint_{t^-_{1,2}\le t_{k+1}\le\cdots\le t_n}
e^{2i[\phi(t_n,\zeta)-\phi(t_{n-1},\zeta)+\cdots 
+\phi(t_{k+1},\zeta)-\phi(t^-_{1,2},\zeta)]}
\prod_{j=k+1}^n f(t_j)\,dt_j
\Big)
\\
\cdot
\Big(
\iint_{x\le t_1\le\cdots\le t_k\le t^+_{1,1} }
e^{2i[\phi(t^+_{1,1},\zeta)
-\phi(t_k,\zeta) + \phi(t_{k-1},\zeta) + \cdots
+ (-1)^{n-1}\phi(t_1,\zeta)]}
\prod_{j=1}^k f(t_j)\,dt_j
\Big).
\end{multline*}

Suppose now that $n$ is even.
The proof in \cite{christkiselevfiltrations}
is an induction based on a repeated application of this decomposition;
each step of that recursion involves a ``cut point'' $t^+_{m,j}
= t^-_{m,j+1}$ playing the
same role as $t^+_{1,1}=t^-_{1,2}$ in the above formula.
At each step, the region of integration is decomposed into 
subregions as above, and corresponding to each subregion there is
a splitting of the terms in the phase into two subsets.
At any step which results in an odd number of terms appearing
in one (hence both) subsets,
we modify the resulting phases by introducing a factor
$1 = \exp(\pm 2i[\phi(t^+_{m,j},\zeta)-\phi(t^-_{m,j+1},\zeta)])$,
factoring it as a product of one exponentials, and splitting those
two exponential factors as above. This, together with the argument
in \cite{christkiselevfiltrations}, yields the assertion of the lemma
for even $n$.

For odd $n$ we introduce a factor of
$1 = \exp(2i[\phi(x,\zeta)-\phi(x,\zeta)])$,
incorporate $\exp(-2i\phi(x,\zeta))$ into the phase, thus
reducing matters again to the case where there an even number of
terms. The remaining factor of $\exp(2i\phi(x,\zeta))$
is bounded above, uniformly for $0\le\Im(\zeta)\le 1$ and $x\ge 0$,
so is harmless. 
\end{proof}

\begin{proof}[Proof of Corollary~\ref{cor:nontangentialmaxbound}]
Let $\Lambda\Subset(0,\infty)$ be a compact interval, and let $1< p<2$,
the case $p=1$ being trivial.
Set $q = p'=p/(p-1)$.
We discuss only the contributions of terms involving 
$\phi(t^-_{m,j},\zeta)$
to $G$ and $G_M$; those involving $t^+_{m,j}$ are treated in 
exactly the same way.
For any $m\ge 1$ and any $f\in L^p$ we have, since $q/2\ge 1$, 
\begin{multline}
\|\Big(
\sum_{j=1}^{2^m} 
\Big|\int_{E^m_j} e^{2i[\phi(t,\lambda+i\e)-\phi(t^-_{m,j},\lambda+i\e)]}
f(t)\,dt\Big|^2
\Big)^{1/2}\|_{L^{q}(\Lambda,d\lambda)}^q
\\
 \le
\Big(\sum_j 
\Big[\int_{\Lambda}
\Big|\int_{E^m_j} e^{2i[\phi(t,\lambda+i\e)-\phi(t^-_{m,j},\lambda+i\e)]}
f(t)\,dt\Big|^q
\,d\lambda\Big]^{2/q}\Big)^{q/2},
\end{multline}
by Minkowski's integral inequality.
By Lemma~\ref{lemma:parseval}, the right-hand side is
\begin{equation*}
\le 
C \Big(\sum_j \|f\cdot \chi_{E^m_j}\|_{L^p}^2\Big)^{q/2}.
\end{equation*}
Since $p\le 2$, this is
\begin{equation*}
\le
C\Big(\sum_j \|f\|_{L^p}^{2-p} \|f\cdot \chi_{E^m_j}\|_{L^p}^p )^{q/2}
= C\|f\|_{L^p}^q. 
\end{equation*}
If we assume that the martingale structure is adapted to $f$ in $L^p$, 
then for $p<2$ we have the improved majorization
\begin{equation*}
\|f\cdot \chi_{E^m_j}\|_{L^p}^2
\le 2^{-m(2-p)/p}\|f\|_{L^p}^{2-p}\|f\cdot \chi_{E^m_j}\|_{L^p}^p,
\end{equation*} 
whence the final bound is $2^{-\rho m}\|f\|_{L^p}^q$ for some
$\rho(p)>0$.

Therefore 
\begin{equation}
\int_\Lambda \|\frakG_M(f)(\lambda+i\e)\|_{\scriptb}^q\,d\lambda
\le C2^{-\rho M}\|f\|_{L^p}^q,
\end{equation}
uniformly for all $\e>0$.
The first two conclusions of the Corollary now follow from  
Lemma~\ref{analyticf},
since $\Lambda$ is an arbitrary compact interval.

That $\frakG(f)(\zeta)$ converges almost everywhere to 
$\frakG(f)(\lambda)$
in the $\scriptb$ norm as $\zeta\to \lambda$ nontangentially, is an 
immediate consequence of the bound 
$\|NG_M(f)\|_{L^q}\le C2^{-\varepsilon M}\|f\|_{L^p}$, since  
\begin{equation}
\int_{t\ge t'} e^{2i[\phi(t, \zeta) - \phi(t', \zeta)]}
f(t)\,dt
\to 
\int_{t\ge t'} e^{2i[\phi(t, \lambda) - \phi(t', \lambda)]}
f(t)\,dt
\end{equation}
almost everwhere as $\zeta\to\lambda$ nontangentially, for all $f\in L^p$.
This holds for all $f$ in the dense subspace $L^1\cap L^p$,
and then follows for general $f$ by Lemma~\ref{analyticf},
Lemma~\ref{lemma:parseval}, and standard reasoning.
\end{proof}

\section{Resolvents and spectral projections}
\label{section:limiting absorption}

\subsection{The half-line case}  \label{section:resolvents}

Consider the operator
\begin{equation}\label{so}
H_V = -\frac{d^2}{dx^2} +V(x)
\end{equation}
on $\reals^+$, with Dirichlet boundary condition at the origin.
Any other selfadjoint boundary condition can be treated in a similar way. 
For each $z\in\complex$, let $u_1(x,z), u_2(x,z)$ be the
unique solutions of
\begin{equation}\label{efeq}
-u''+V(x)u = zu
\end{equation}
satisfying the boundary conditions
$u_1(0,z) = (0,1)^{t}, u_2(0,z)=(1,0)^t$;
the superscripts $t$ denote transposes.
The classical theory of second order differential operators
(see e.g. \cite{CL,Tit}) tells us that if \eqref{so} is in the
limit point case, then for any $z \in \complex \setminus \reals$
there exists a unique complex number
$m(z)$, called the Weyl $m$-function, such that
\[ f(x,z) = u_1(x,z)m(z)+u_2(x,z) \in L^2(\reals^+). \]
We will always consider potentials which
lead to the limit point situation,
as is the case if $V \in L^1+L^p$ for some $1\le p<\infty$; see for example
\cite{ReSi2}.
The Weyl $m$ function is a Herglotz function, that is, it is analytic
in the upper half-plane and has positive
imaginary part there. 

By direct computation, the resolvent $(H_V-z)^{-1}g$
for $z\in\complex^+$ is given by
\begin{equation}\label{resol}
(H_V-z)^{-1} g(x) = u_1(x,z) \int_x^\infty f(y,z)g(y)\,dy +
f(x,z) \int_0^x u_1(y,z) g(y)\,dy.
\end{equation}
Denote by $P_{(a,b)}$ the spectral projection associated to
$H_V$ and to the interval $(a,b)$.
>From the resolvent formula \eqref{resol} we may derive a
formula for the projections $P_{(a,b)}$. In doing so, we will use the
following three routine facts.
\\
1. The functions $u_1(x,z), u_2(x,z)$ are continuous in $x$ for
each $z$, and are entire holomorphic functions of $z$ for each $x$.
\\
2. The $m$ function (in fact, any Herglotz function) 
has a representation
\begin{equation}
m(z) = C_1 + C_2z+
\int_\reals \left( \frac{1}{t-z} - \frac{t}{1+t^2} \right)
d\mu(t)
\end{equation}
for some positive Borel measure $\mu$ satisfying $\int (1+|t|^2)^{-1}
d\mu(t) < \infty,$  see e.g. \cite{AD}.
In the $m$ function context, $\mu$ is often called the
spectral measure. For Dirichlet boundary conditions, as considered here,
it is the spectral measure corresponding to the generalized vector $\delta_0'$,
defined by $\delta_0'(u)=u'(0)$ for any function $u$ in the domain of $H_V.$  
The moment condition above corresponds to the fact that
the derivative $\delta'_0$ belongs to the Sobolev-like space $H_{-2}(H_V)$ associated
to $H_V$ (see e.g.\ \cite{BeSh} for details on families of Sobolev-like spaces 
associated with any selfadjoint operator $A$).  \\
3. $\Im m(E+i\epsilon)$ converges weakly to $\pi \mu$
as $\epsilon \rightarrow
0^+$. Moreover, $\Im m(E+i\epsilon)$ has limiting boundary values for
Lebesgue-almost every $E$,
and the density of the absolutely continuous part of $\mu$ satisfies
\begin{equation}\label{acpart}
d\mu_{\rm{ac}}(E)= \frac{1}{\pi}\Im m(E+i0) \,dE,
\end{equation}
where $m(E+i0) = \lim_{\e\to 0^+}m(E+i\e)$.
Since the imaginary part of $m$ is simply a Poisson integral of $\mu$,
this is straightforward.

Fix functions $h$ and $g$ with compact support.
We integrate the resolvent element $\langle (H_V-z)^{-1}g,h \rangle$
over the contour $\gamma_\epsilon$ in complex plane consisting of two
horizontal intervals $(a \pm i \epsilon, b \pm i \epsilon)$, and
two vertical intervals at the ends
connecting them. In the limit $\epsilon \rightarrow 0^+$,
the contributions of the vertical intervals disappear unless $a$ or $b$
is an eigenvalue (point mass of $\mu$); we will assume this
is not the case.
By the spectral theory (see, e.g. \cite{ReSi1}, Stone formula)
we get then the following
expression for $\langle P_{(a,b)}g,h \rangle$:
\begin{multline}\label{spro}
\langle P_{(a,b)}g,h \rangle
= \frac{1}{2\pi i}\lim_{\epsilon \rightarrow 0}
\int_{\gamma_\epsilon}
 \langle (H_V-z)^{-1}g,h \rangle \, dz  \\
=
\int_a^b  \int_\reals   \int_\reals
u_1(x,E)u_1(y,E) g(x)\overline{h}(y)\,dx\,dy
\,d\mu(E).
\end{multline}
In passing to the last line,
we have taken into account the resolvent formula \eqref{resol},
the properties of $u_{1},u_2$
(in particular, $u_2$ drops out because of analyticity), and the fact that
$\pi^{-1}\Im m(E+i\epsilon)$ converges weakly to $\mu$. The compact
supports of $h,g$ ensure that the integral is well-defined.
Similar formulas can be found in \cite{CL, Tit}.

Since $P_{(a,b)}$ is by its definition an orthogonal projection,
an immediate consequence
of \eqref{spro} is that the mapping $g\mapsto \int_\reals
u_1(x,E)g(x)\,dx$, initially defined for continuous $g$ having compact
support, extends to an orthogonal projection from $L^2(\reals^+,dx)$
to $L^2((a,b),\mu)$.
Dually, from \eqref{spro} we see that for each $g\in L^2(\reals,d\mu)$,
\begin{equation}
U_0 g(x) =\lim_{N \rightarrow \infty} \int_{-N}^N
u_1(x,E)g(E) d\mu(E)
\end{equation}
exists in $L^2(\reals^+,dx)$ norm,
and that the linear operator $U_0$ thus defined
is a unitary bijection from $L^2(\reals,d\mu)$ to $L^2(\reals,dx)$
with inverse
\begin{equation}
 U_0^{-1}
g(E) = \lim_{N \rightarrow \infty} \int_{-N}^N u_1(x,E) g(x) dx,
\end{equation}
where the limit is again taken in $L^2$ norm.

With the formula \eqref{spro} for the spectral projections in hand,
we can invoke general spectral theory to find expressions for
the spectral representation and other functions of $H_V$
(see, e.g.\ \cite{BS}).
To $H_V$ and any interval $(a,b)$
are associated a maximal closed subspace of $L^2(\reals^+)$
on which $H_V$ has purely absolutely continuous spectrum,
and spectrum contained in $(a,b)$.
By \eqref{spro} and \eqref{acpart} as well as by
definition of the absolutely continuous part of the spectral projection,
the projection $P^{\rm{ac}}_{(a,b)}$ of $L^2(\reals^+)$
onto this subspace can be written as
\begin{equation}\label{acpro}
P_{(a,b)}^{\rm{ac}}g(x) =
\frac{1}{\pi}
\int_a^b u_1(x,E)  \Big( \int_\reals u_1(y,E) g(y) \,dy \Big)
\Im m(E+i0)\,dE.
\end{equation}
The integral over $\reals$ here is generally understood in the 
$L^2$-limiting sense; we will omit such explanatory remarks in the future.
Consequently the
operator $U$ mapping continuous functions with compact support to
$L^2(\reals, dx)$, defined by
\begin{equation} Uh(x) =  \frac{1}{\pi}\int_{\reals}
u_1(x,E) h(E) \Im m(E+i0) \,dE \end{equation}
extends to an isometry of $L^2(\reals,\Im m(E+i0) \,dE)$
onto the absolutely continuous subspace associated to $H_V$.

The unitary evolution operator on the absolutely continuous subspace is
given by
\begin{equation}\label{evol}
e^{-iH_V t} g(x) =
\frac{1}{\pi}\int_\reals e^{-iEt}u_1(x,E) \tilde{g}(E)
\Im m(E+i0)\,dE,
\end{equation}
where
\[ \tilde{g}(E) =
\frac{1}{\pi}\int_\reals u_1(y,E) g(y)\,dy. \]
Finally, in the case $V=0$ we can compute explicitly
$u_1(x,E) = \sqrt{E}^{-1} \sin \sqrt{E}x,$
$m(z) = \sqrt{z},$ and so the evolution operator can be written as
\begin{equation}\label{freeevol}
e^{-iH_0 t} g(x) =
\frac{1}{\pi}\int_\reals e^{-iEt} \sin(\sqrt{E}x) \hat{g}(E)
\,\frac{dE}{\sqrt{E}},
\end{equation}
where
\[ \hat{g}(E) = \int \sin (\sqrt{E}x) g(x)\,dx. \]

For almost every $E>0$, define the scattering coefficient
$\gamma(E)\in\complex$ by
\begin{equation} \label{reflectioncoefficientdefn}
\gamma(E) = 1/ u(0,E)
\end{equation}
where $u(x,E)$ is the unique generalized eigenfunction
asymptotic to $\exp(i\xi(x,E))$ as $x\to+\infty$,
whose existence was established in Theorem~\ref{thm:complexeigenfunctions}.

The following proposition connects the formulae of this section
with the generalized eigenfunctions analyzed in \S\ref{section:complex}.
\begin{proposition}[Limiting absorption principle]\label{relim}
Assume that $V\in L^1+L^p(\reals)$ for some $1<p<2$.
For almost every $E \in \reals^+$, the generalized eigenfunction
$f(x,E+i0)=u_1(x,E)m(E+i0)+ u_2(x,E)$ satisfies
\begin{equation}\label{asbes}
f(x,E+i0) =  \gamma(E) e^{i\xi(x,E)}(1+o(1)),
\end{equation}
and
\begin{equation} |\gamma(E)|^2 = \Im m(E+i0)/\sqrt{E}.
\end{equation}
\end{proposition}

\begin{proof}[Proof of Proposition~\ref{relim}]
Denote by $u(x,z)$ the generalized eigenfunctions, with spectral parameter $z$,
whose existence was established in Theorem~\ref{thm:complexeigenfunctions}.
By the uniqueness of $L^2$ solutions of \eqref{eigenfneqn}
for $z \in \complex^+$,
$f(x,z) = u(x,z)/u(0,z).$ For a.e.~$E,$ $f(x,z)$ converges
uniformly as a function of $x\in\reals^+$ to
$u_1(x,E)m(E+i0)+u_2(x,E)$, as $z$ converges to $E$ nontangentially.
The uniform convergence ensures that $f(x,E+i0)$ is indeed
a generalized eigenfunction associated to the spectral parameter $E$.
At the same time,
\[ \frac{u(x,E+i\epsilon)}{u(0,E+i\epsilon)} \rightarrow
\frac{u(x,E)}{u(0,E)} \]
for almost every $E$ by Theorem~\ref{thm:complexeigenfunctions},
and therefore for $E>0$
\[ u_1(x,E)m(E+i0)+u_2(x,E)=\frac{u(x,E)}{u(0,E)}= \gamma(E)
e^{i \xi(x,E)}(1+o(1)) \]
as $x \rightarrow \infty$
(there is no absolutely continuous spectrum for $E<0$).
The relation between $\gamma$ and $m(E+i0)$
follows by comparing the Wronskians of the left
and right hand sides (taken with their complex conjugates).
\end{proof}


Now we are going to rewrite the spectral representation in a manner
convenient for the scattering theory. Notice that
\[ u_1(x,E) = \frac{1}{2i \Im m(E+i0)}
\left( \gamma(E)u(x,E) - \ov{\gamma}(E)\ov{u}(x,E)
\right). \]
Introduce
\begin{equation}
\psi(x,\lambda) = \lambda \ov{\gamma}(\lambda^2) u_1(x,\lambda^2)
= \frac{1}{2i}\Big(
u(x,\lambda^2)
- \frac{{\ov{\gamma}(\lambda^2)}}{\gamma(\lambda^2)}
\cdot \ov{u}(x,\lambda^2) \Big).
\end{equation}
where $u(x,E)$ continue to denote the generalized eigenfunctions
whose existence was established in Theorem~\ref{thm:complexeigenfunctions},
now for $E\in\reals^+$.
Then the results of this subsection may be summarized as follows.
\begin{proposition}
Suppose that $V\in L^1+L^p(\reals^+)$ for some $1<p<2$.
Then for the associated selfadjoint Schr\"odinger operator $H_V$
on $L^2(\reals^+)$ with Dirichlet boundary conditions,
the spectral projection $P_{(a,b)}^{{\rm ac}}$ can be expressed as
\begin{equation}
P_{(a^2,b^2)}^{{\rm ac}}g(x)= \frac{2}{\pi}
\int_{[a,b]\cap\reals^+} \psi(x,\lambda)
\tilde g(\lambda)
\,d\lambda
\end{equation}
for any $g\in L^2(\reals^+)$, where
the modified Fourier transform $\tilde g$ is defined by
\begin{equation} \tilde{g}(\lambda) = \int_\reals \ov{\psi}(x,\lambda)
g(x)\,dx. \end{equation}
Similarly, the associated wave group is
\begin{equation}\label{perev}
e^{-it H_V} g(x) =
\frac{2}{\pi}\int\limits_0^\infty e^{-i\lambda t}
\psi(x,\lambda)\tilde g(\lambda) \,d\lambda.
\end{equation}
The operator
\begin{equation}
U_V f(x) = \sqrt{2/\pi} \int_0^\infty f(\lambda) \psi(x,\lambda)\,d\lambda
\end{equation}
is a unitary surjection from $L^2(\reals^+,d\lambda)$
onto $\scripth_{\text{ac}}$.
\end{proposition}


\subsection{Formulae for the case of the whole real line}


Consider $-d^2/dx^2+V(x)$ as an essentially selfadjoint operator
on $L^2(\reals)$.
The asymptotic analysis of the preceding section, 
in which $x\to+\infty$, works equally well as $x\to-\infty$.
Whereas it was necessary in the half-line case to 
relate asymptotics at $+\infty$ to boundary conditions at $x=0$,
now we must relate asymptotics at $+\infty$ to those at $-\infty$.
The absolutely continuous spectrum now has multiplicity two,
which complicates some of the formulae.

Let $z\in\complex^+$. Solutions
$u_1(x,z),u_2(x,z)$ are defined precisely as before, with the
same initial conditions at $x=0$; now they
are considered as global solutions on $\reals$ rather than merely
on $[0,\infty)$.
Introduce two solutions
$f_\pm(x,z) = u_1(x,z) m_\pm(z) + u_2(x,z)$, 
so that $f_+\in L^2(\reals^+,dx)$ and
$f_-\in L^2(\reals^-,dx)$ for each $z\in\complex^+$.
In terms of these, the resolvent is given by
\begin{equation} 
(H_V-z)^{-1}g(x) = \frac{-1}{W[f_+,f_-])}
\Big( f_+(x,z) \int\limits_{-\infty}^x f_-(y,z) g(y) \,dy +
f_-(x,z) \int\limits_{x}^{\infty} f_+(y,z)g(y)\,dy \Big). 
\end{equation}
Notice that $W[f_+,f_-]=m_+-m_-$ 
(since $\Im m_-(z)<0$ in $\complex^+$).

The formula for the spectral projection associated to the absolutely
continuous spectrum can be computed from the resolvent in a way similar 
to the half-line case.

In the free case $V=0,$ it simplifies to
\begin{equation}\label{freewa}
P_{(a^2,b^2)}^{{\rm ac}}g(y) = \frac{1}{2\pi}
\int\limits_a^b \chi_{\reals^+}(\lambda) \bigg(
e^{i\lambda x} \int\limits_\reals e^{-i\lambda y }g(y)\,dy +
e^{-i\lambda x} \int\limits_\reals e^{i\lambda y }g(y)\,dy \bigg)
\ d\lambda .
\end{equation}
The evolution operator is given by
\begin{equation}\label{frevwa}
e^{-iH_0t}g(x) = \frac{1}{2\pi}
\int\limits_0^\infty  e^{-i\lambda^2 t}\left(
e^{i\lambda x} g_+(\lambda) +
e^{-i\lambda x} g_-(\lambda) \right)\, d\lambda,
\end{equation}
where
\[ g_\pm (\lambda) = \int\limits_\reals e^{\mp i \lambda y} g(y) dy. \]

In the general case, we introduce scattered waves,
defining $\psi_\pm(x,\lambda)$ for almost every $\lambda \in\reals^+$
to be the unique generalized eigenfunctions with the asymptotic
behavior
\begin{equation}\label{risc}
\psi_+(x,\lambda) = \left\{ \begin{array}{ll} t_1(\lambda) e^{i\phi(x,\lambda ) } (1+o(1)), 
& x \rightarrow +\infty, \\
\left( e^{i\phi(x,\lambda )} +r_1(\lambda )e^{-i\phi(x,\lambda )} \right)\ +o(1), 
& x \rightarrow -\infty \end{array} \right.
\end{equation}
\begin{equation}\label{lesc}
\psi_-(x,\lambda) = \left\{ \begin{array}{ll} t_2(\lambda) e^{-i\phi(x,\lambda)} (1+o(1)), 
& x \rightarrow -\infty, \\
\left( e^{-i\phi(x,\lambda)} +r_2(\lambda)e^{i\phi(x,\lambda)} \right)\ +o(1), 
& x \rightarrow +\infty, \end{array} \right.
\end{equation}
where we recall that $\phi(x,\lambda)= \lambda x - (2\lambda)^{-1} \int_0^x V(s)\,ds.$  
That is, $\psi_\pm$ are defined by the stated asymptotics as
$x\to\pm\infty$, respectively;
then writing them as linear combinations of the generalized
eigenfunctions $\exp(\pm i\xi(x,\lambda^2))+o(1)$ as $x\to\mp\infty$,
respectively, we define $r_i,t_i$ to be the scalar coefficients
in those linear combinations.
The asymptotics in \eqref{risc}, \eqref{lesc} can be differentiated,
in the sense that $\psi_+(x,\lambda) = i\lambda t_1(\lambda)\exp(i\xi(x,\lambda^2))
+ o(1)$ as $x\to+\infty$, with analogous formulae as $x\to-\infty$
and for $\psi_-$.

As in the half-line case, there exist coefficients $\gamma_\pm(\lambda)$
such that
\begin{equation}\label{rel1}
\psi_\pm (x,\lambda) = \gamma_\pm^{-1}(\lambda) (u_1(x,\lambda^2) 
m_\pm(\lambda^2+i0)+u_2(x,\lambda^2))
\end{equation}
for almost every $\lambda$. 
Computations of Wronskians,  
examination of transfer matrices
between $\psi_+,\ov{\psi}_+$ and $\psi_-,\ov{\psi}_-$,
and exploitation of complex conjugation
leads to the following 
formulae relating the transmission and reflection
coefficients $t_i,r_i$ with one another,
and with the coefficients $\gamma_\pm,m_\pm$:
\begin{align}
&|r_i|^2+|t_i|^2=1 \text{ for } i=1,2,
\label{traref}
\\
&t_1=t_2, 
\\
&r_2=-\frac{t_1}{\ov{t}_1}\ov{r}_1,
\\
&|\gamma_\pm|^2|t_1|^2 \lambda = \pm\Im m_\pm(\lambda^2+i0). 
\label{gamm}
\end{align}
%
After a computation, we obtain
\begin{equation}\label{spscwa}
P_{(a^2,b^2)}^{{\rm ac}}g(y) =  \frac{1}{2\pi}
\int\limits_a^b \chi_{\reals^+}(\lambda) \bigg(
\psi_+(x,\lambda) \int\limits_\reals \ov{\psi}_+(y,\lambda ) g(y)\,dy +
\psi_-(x,\lambda ) \int\limits_\reals \ov{\psi}_-(y,\lambda ) g(y)\,dy \bigg)
\ d\lambda .
\end{equation}
Therefore, the evolution operator is given by
\begin{equation}\label{perevwa}
e^{-iH_V t}g(x) = \frac{1}{2\pi}
\int\limits_0^\infty  e^{-i \lambda^2 t}\Big(
\psi_+(x,\lambda ) \tilde{g}_+(\lambda ) +
\psi_-(x,\lambda ) \tilde{g}_-(\lambda ) \Big)\ d\lambda ,
\end{equation}
where the transforms $\tilde g_\pm$ are defined by
\begin{equation}
\tilde{g}_\pm (\lambda) = \int\limits_\reals \ov{\psi}_{\pm}(y,\lambda) g(y) dy. 
\end{equation}

\subsection{Dirac-type operators on the whole real line} 

Here we consider a
system of differential equations of second order which plays an important 
role in analysis of the defocusing
nonlinear Schr\"odinger equation (NLS) by the inverse scattering method.
This system is given by \cite{Nov}
\begin{equation}\label{orsys}
y' =  iz 
\begin{pmatrix}
1 & 0 \\ 0 & -1 
\end{pmatrix} y
+ 
\begin{pmatrix}
0 & q \\ \ov{q} & 0 
\end{pmatrix} y.
\end{equation}

We are going to use two alternative equivalent representations of this 
system.  The first is
\begin{equation}\label{forsol}
\begin{pmatrix}
-i\partial_x & V \\ \ov{V} & i \partial_x 
\end{pmatrix}
y  = zy,
\end{equation}
where $V = iq$. System \eqref{forsol} is obtained from \eqref{orsys} by 
multiplying both sides by a matrix 
$\left( \begin{array}{cc} 1 & 0 \\ 0 & -1 \end{array} \right)$.
The second is obtained by setting $y=Q^{-1}\phi$ in \eqref{forsol}, 
where
\[ Q =\left( \begin{array}{cc} 1 & 1 \\ i & -i \end{array} \right). \]
We obtain
\begin{equation}\label{dirac}
\left( \begin{array}{cc} 0 & -\partial_x \\ \partial_x & 0 \end{array} 
\right) 
\phi +
\left( \begin{array}{cc} \Re V & \Im V \\ \Im V & -\Re V \end{array} 
\right) \phi = z \phi.
\end{equation}
The operator on the left hand side of \eqref{dirac} is a particular case 
of a second order Dirac-type operator (the most general case does not 
require the second diagonal entry to be minus the first; in the most 
general case there will be a phase shift in the main term of the solution
asymptotics, unlike the situation here).

We digress briefly to comment on  the 
possibility of embedded point spectrum for this operator. 
Let $z=E \in \reals.$ The equation \eqref{forsol} preserves the analogue 
of the Wronskian for two solutions $f=(f_1,f_2)$ and $g=(g_1,g_2),$ 
$W[f,g]=i(f_2g_1-f_1g_2).$
Notice that if $(g_1, g_2)$ is a solution, so is $(\ov{g}_2, \ov{g}_1).$ 
The constancy of the Wronskian for these two solutions implies that 
$|g_2|^2-|g_1|^2$ is constant. Consider for simplicity the
case where $|g_1|=|g_2|$.
Representing $g_1=R(x) e^{i\theta_1(x)},$ $g_2=R(x) e^{i\theta_2(x)}$ 
and using \eqref{forsol}, we find that
$\theta_1+\theta_2=c$ is constant, and
\begin{equation}\begin{split}  \label{pruf}
(\log R)' &= - \Re V \sin(2\theta_1 -c)+ \Im V \cos(2\theta_1-c) \\
\theta'_1 &= E - \Re V \cos (2\theta_1-c) -\Im V \sin (2\theta_1 -c).
\end{split}\end{equation}
{}From these Pr\"ufer-like equations one can easily recover 
certain constructions of embedded eigenvalues available in the
Schr\"odinger case \cite{Nab,KLS,KRS}; 
see \cite{Nab} for an earlier alternative approach to the Dirac case.
In particular, there exist
potentials $V=O(|x|\rp)$ with isolated eigenvalues embedded in
the continuous spectrum, and for any function $g(|x|)$ tending
to infinity, there exists $|V(x)|\le g(|x|)/|x|$
for which there is a dense set of eigenvalues in $\reals^+$.

The form \eqref{forsol} will prove useful for studying solutions; 
the form \eqref{dirac} will
allow us to set up scattering in a way completely parallel to 
the setup for Schr\"odinger operators on the whole axis.
First, setting
\[ y_1 = \left( \begin{array}{cc} e^{izx} & 0 \\ 0 & e^{-izx} 
\end{array} \right) y \]
we arrive at
\[ y_1' = \left( 
\begin{array}{cc} 0 & -iV e^{-2izx} 
\\ 
i\ov{V} e^{2izx} & 0 \end{array} \right) y_1. \]
Iterating as in the Schr\"odinger case, we get  
solutions $\psi_\pm(x,z)$ for any  $z \in \complex^+$ and 
almost every $z \in \reals$ with the asymptotic behavior
\begin{eqnarray}\label{psias}
\psi_+(x,z) = e^{izx} 
\left( \left( \begin{array}{c} 1 \\ 0 \end{array} \right)
+ o(1) \right), & x \rightarrow +\infty \\
\nonumber
\psi_-(x,z) = e^{-izx} 
\left( \left( \begin{array}{c} 0 \\ 1 \end{array} \right)
 + o(1) \right), & x \rightarrow -\infty.
\end{eqnarray}
An analogue of Proposition~\ref{relim} also holds.

Before proceeding, we record a corollary of the existence of such generalized eigenfunctions
for real $z$.
\begin{corollary}  \label{cor:diracspectrum}
For any $1<p<2$ and any $V\in L^1+L^p(\reals)$, the Dirac operator $D_V$
on $L^2(\reals)$ has nonempty absolutely continuous spectrum. More precisely, an essential 
support for the ac spectrum is the entire real line $\reals$. 
\end{corollary}
As remarked above, this hypothesis does not preclude the presence of a dense set of eigenvalues
embedded in the continuous spectrum.

Now let us return to \eqref{dirac}. 
The existence of solution \eqref{psias} and the relation between 
\eqref{dirac} and \eqref{forsol} imply existence of solutions 
$\eta_\pm(x,z)$ of \eqref{dirac},
for $z \in \complex^+ \cup \reals$, of the form
\begin{eqnarray}\label{phias}
\eta_+(x,z) = e^{izx} 
\left( \left( \begin{array}{c} 1 \\ i \end{array} \right)
+ o(1) \right), & x \rightarrow +\infty \\
\nonumber
\eta_-(x,z) = e^{-izx} 
\left( \left( \begin{array}{c} 1 \\ -i \end{array} \right) + o(1) \right), & x \rightarrow -\infty.
\end{eqnarray}
We proceed very much as in the Schr\"odinger case.
The role of the Wronskian is played by $W[f,g]=f_2g_1-f_1g_2$.
As before, we denote by $u_{1,2}(x,z)$ solutions satisfying
$u_1(0,z)=(0,1)^T,$ $u_2(0,z)=(1,0)^T.$
Let $f_{\pm}(x,z)=u_1(x,z)m_{\pm}(z)+u_2(x,z)$ be $L^2$ solutions
on $\pm \infty.$
Denote the operator on the left hand side of
\eqref{dirac} by $D_V.$ Its resolvent is given by
\begin{equation}\label{dvres}
(D_V - z)^{-1} \left( \begin{array}{c} g_1 \\ g_2 \end{array} \right)
= \frac{1}{W[f_-,f_+]} \left( f_+ \int\limits_{-\infty}^x
(f_{-,1} g_1+f_{-,2}g_2) \, dy + f_-
\int\limits_x^\infty  (f_{+,1}g_1+f_{+,2}g_2) \,dy
\right).
\end{equation}
Let us denote in this section
\[ \langle f(x), g(x) \rangle = \int\limits_\reals
(f_1(x)g_1(x)+f_2(x)g_2(x))\,dx;  \]
when necessary, the integral is understood in a limiting sense
(in $L^2$ if $f$ is a solution which also depends on $E$).
An analog of \eqref{acpro}
 can be derived using contour integration.
The conclusion is that (one version of)
the spectral representation of the absolutely continuous part of the
spectral measure is given by the map
\[  
U: L^2(\complex^2,dx) \mapsto L^2(\complex^2, MdE), \,\,\,
U \left( \begin{array}{c} g_1 \\ g_2 \end{array} \right) =
 \left( \begin{array}{c} \tilde{g}_1 \\ \tilde{g}_2 \end{array} \right), 
\]
 where
\[ \tilde{g}_i (E) =
\langle u_i(x,E), g \rangle  \]
and
\[ M(E) = \Im \left( \begin{array}{cc}  \frac{m_+ m_-}{m_- - m_+} &
 \frac{m_-}{m_- - m_+} \\  \frac{m_-}{m_--m_+} & \frac{1}{m_--m_+}
\end{array} \right)(E+i0). \]
Then $U$ is an isometry, $U D_V U^{-1} =E$
(with $U^{-1}$ given by
\[ U^{-1} \tilde{g} (x) =
\int\limits_\reals 
\begin{pmatrix} u_1(x,E)\\ u_2(x,E) \end{pmatrix}
M \tilde{g} \,dE, \]
where $(u_1,u_2)$ is a $2 \times 2$ matrix with columns $u_1,$ $u_2.$
The unitary evolution of the absolutely continuous part is then given by
\begin{equation}\label{dev}
e^{-iD_V t} g(x) =  \int\limits_\reals
e^{-iEt} 
\begin{pmatrix} u_1(x,E)\\ u_2(x,E) \end{pmatrix}
M \tilde{g}(E) \,dE.
\end{equation}
In the free case, the spectral representation can be simplified.
In particular, we get
\begin{equation}\label{dfrev}
e^{-iD_0t}g(x) = \frac{1}{4} \int\limits e^{-iEt} \left(
\eta_{+,0}(x,E) \langle \ov{\eta}_{+,0}, g \rangle +
\eta_{-,0}(x,E) \langle \ov{\eta}_{-,0}, g \rangle \right),
\end{equation}
where
\[ \eta_{\pm,0}(x,E) = e^{\pm iEx} \left( \begin{array}{c}
1 \\ \pm i \end{array} \right). \]
In the perturbed case, we introduce the scattered waves
\begin{equation}\label{drisc}
\eta_+(x,E) = \left\{
\begin{array}{ll} t_1(E) e^{iEx}
 \left( \begin{array}{c}
1 \\  i \end{array} \right) +o(1), & x \rightarrow +\infty, \\
\left( e^{iEx} \left( \begin{array}{c}
1 \\  i \end{array} \right) +r_1(E)e^{-iEx} \left( \begin{array}{c}
1 \\ - i \end{array} \right) \right)+o(1), & x \rightarrow -\infty
 \end{array} \right.
\end{equation}
\begin{equation}\label{dlesc}
\eta_-(x,E) = \left\{ \begin{array}{ll} t_2(E) e^{-iEx}
 \left( \begin{array}{c}
1 \\ - i \end{array} \right) +o(1), & x \rightarrow -\infty, \\
\left( e^{-iEx} \left( \begin{array}{c}
1 \\ - i \end{array} \right)
 +r_2(E)e^{iEx} \left( \begin{array}{c}
1 \\  i \end{array} \right)
 \right)+o(1), & x \rightarrow +\infty \end{array} \right.
\end{equation}
Existence of such solutions follows from \eqref{psias}.
A computation parallel to the whole axis Schr\"odinger 
operator case allows us to rewrite the dynamics \eqref{dev} as
\begin{equation}\label{dscev}
e^{-iD_V t}g(x) = \frac{1}{4} \int\limits_\reals e^{-iEt} \left(
\eta_{+}(x,E) \langle \ov{\eta}_{+}, g \rangle +
\eta_{-}(x,E) \langle \ov{\eta}_{-}, g \rangle \right).
\end{equation}

\section{Long-time asymptotics} \label{section:heart}

Let $V\in L^1+L^p(\reals^+)$ for some $1<p<2$,
and consider the selfadjoint 
Schr\"odinger operator $H_V=-d^2/dx^2+V(x)$
on $L^2(\reals^+)$ with Dirichlet boundary condition.
Let $\scripth=L^2(\reals^+)$ and let
$\scripth_{\text{ac}}$ denote the maximal closed subspace
of $\scripth$ on which $H_V$ has purely absolutely continuous
spectrum.

Let $\phi(x,\lambda) = \lambda x-(2\lambda)\rp\int_0^x V$
be as before,
and for almost every $\lambda\in\reals^+$ let $u(x,\lambda)$
be the unique generalized eigenfunction associated to
the spectral parameter $\lambda^2$, satisfying 
 $\Psi(x,\lambda) = \exp(i\phi(x,\lambda))+o(1)$ as $x\to+\infty$.
>From the preceding section recall the generalized eigenfunctions
$\psi(x,\lambda) = (2i)\rp(u - e^{i\omega(\lambda)}\ov{u}))$,
where $\exp(i\omega(\lambda)) = \ov\gamma(\lambda^2)/\gamma(\lambda^2)$.
Recall the unitary bijection $U_V:L^2(\reals^+,d\lambda)\mapsto
\scripth_{\text{ac}}$ defined
by $U_Vf(x) = \sqrt{2/\pi}\int_0^\infty f(\lambda)\psi(x,\lambda)\,d\lambda$.
Finally, recall that
$e^{-iH_Vt}(U_V f)(x) = \sqrt{2/\pi}\int_0^\infty e^{-i\lambda^2 t}f(\lambda )\psi(x,\lambda)\,
d\lambda$.

Fix a martingale structure $\{E^m_j\}$ on $\reals^+$ that is
adapted to $V$, in the sense that $\int_{E^m_j}|V|^p
 = 2^{-m}\int_{\reals^+}|V|^p$ for all $m\ge 0$ and all $j$.
For any sufficiently small $\delta>0$, recall the functional
\begin{equation*}
g_\delta(f)(\lambda)
= \sum_{m=1}^\infty 2^{m\delta}
\Big(
\sum_{j=1}^{2^m} 
\big|\int_{E^m_j} e^{2i\phi(x,\lambda)}f(x)\,dx\big|^2
+
\big|\int_{E^m_j} e^{-2i\phi(x,\lambda)}f(x)\,dx\big|^2
\Big)^{1/2}.
\end{equation*}
We have shown that for all sufficiently small $\delta_0$,
$g_{\delta_0}(V)(\lambda)<\infty$ for almost every $\lambda\in\reals^+$.
Fix some $\delta<\delta_0$.
\begin{definition}
A compact set $\Lambda\Subset(0,\infty)$
is said to be a {\em set of uniformity} if
$g_\delta(V)$ is a bounded function of $\lambda\in\Lambda$,
and if 
$u(x,\lambda)-e^{i\phi(x,\lambda)}\to 0$ 
as $x\to+\infty$, uniformly for all $\lambda\in\Lambda$.
\end{definition}

\begin{lemma}
For any $f\in L^2(\reals^+)$,
for any $R<\infty$,
\begin{equation*}
\|e^{-itH_V}U_V f\|_{L^2([0,R])} \to 0
\ \  \text{as } t\to\infty.
\end{equation*}
\end{lemma}

\begin{proof}
Since $\exp(-itH_V)$ is unitary, it suffices to prove this 
for all $f$ in some dense subspace of $L^2(\reals^+)$;
we choose the subspace consisting of all $f$ supported on 
some set of uniformity $\Lambda$ (which depends on $f$).
The functions $\psi(x,\lambda)$ are uniformly bounded on
$[0,\infty)\times\Lambda$, as are their derivatives with respect
to $x$, from which it follows that $U_V$ is a compact mapping
from $L^2(\Lambda)$ to $L^2([0,R])$. Thus it suffices to establish
weak convergence to zero. 
Now for any $h\in L^2([0,R])$, 
\begin{equation}
\langle e^{-itH_V}U_V f,h\rangle = \int_\Lambda
e^{-i \lambda^2 t }f(\lambda)U_V^*h(\lambda)\,d\lambda,
\end{equation}
which converges to zero by the Riemann-Lebesgue lemma,
since $U_V^*h\in L^2(\reals^+)$.
\end{proof}

Define
\begin{align}
\psi_0(x,\lambda) & = (2i)\rp\big( e^{i\phi(x,\lambda)}
-e^{i\omega(\lambda)}e^{-i\phi(x,\lambda)}\big)
\\
U_V^\dagger f(x) & = \sqrt{2/\pi} \int_0^\infty f(\lambda) \psi_0(x,\lambda)\,d\lambda;
\end{align}
these are approximations to $\psi,U_V$, respectively.
The next lemma is the key step in showing that only the 
leading-order approximation $\psi_0$ to $\psi$ contributes to the 
long-time asymptotics of the wave group.

\begin{lemma} \label{truekeylemma}
For any set of uniformity $\Lambda$ and any $R>0$,
there exists $C(R,\Lambda)<\infty$ such that for all
$f\in L^\infty(\Lambda)$,
\begin{equation}
\|(U_V-U_V^\dagger)f\|_{L^2([R,\infty))}
\le C(R,\Lambda)\|f\|_{L^\infty},
\ \text{ where }\ C(R,\Lambda)\to 0 \text{ as } R\to\infty.
\end{equation}
\end{lemma}

\begin{proof}
It suffices to prove this with $\psi(x,\lambda)-\psi_0(x,\lambda)$ replaced by
\begin{multline}
\Psi(x,\lambda) - e^{i\phi(x,\lambda)}
\\
= e^{i\phi(x,\lambda)}\sum_{n=1}^\infty S_{2n}(V,V,\dots,V)(x,\lambda)
-  
e^{-i\phi(x,\lambda)}\sum_{n=0}^\infty S_{2n+1}(V,V,\dots,V)(x,\lambda),
\end{multline}
for the conclusion for the other terms follows from this by complex conjugation.
We argue by duality; let $h$ be an arbitrary function in
$L^2([R,\infty))$ of norm $1$, and consider
\begin{multline}
\int_R^\infty e^{i\phi(x,\lambda)}h(x)S_{2n}(V,V,\dots,V)(x,\lambda)
\\
= \iint_{R\le t_0\le t_1\le\cdots\le t_{2n}}
e^{i\phi(t_0,\lambda)}h(t_0)\,dt_0\,
\prod_{k=1}^{2n}e^{\pm_k 2i\phi(t_k,\lambda)} V(t_k)\,dt_k.
\end{multline}
Here $\pm_k$ denotes a plus or minus sign depending on $k$ in any manner;
in fact these signs alternate in our expansion, but that is of no
importance here.
Introduce
\begin{equation}
g_{-\delta}(h)(\lambda)
= \sum_{m=1}^\infty 2^{-m\delta}
\Big(
\sum_{j=1}^{2^m} 
\big|\int_{E^m_j} e^{i\phi(x,\lambda)}h(x)\,dx\big|^2
+
\big|\int_{E^m_j} e^{-i\phi(x,\lambda)}h(x)\,dx\big|^2
\Big)^{1/2}.
\end{equation}
By Proposition~\ref{prop:numericalbound},
\begin{equation}
\Big|\iint_{R\le t_0\le t_1\le\cdots\le t_{2n}}
e^{i\phi(t_0,\lambda)}h(t_0)\,dt_0\,
\prod_{k=1}^{2n} e^{\pm_k 2i\phi(t_k,\lambda)}V(t_k)\,dt_k
\Big|
\le \frac{C^{n+1}}{\sqrt{2n!}}
g_{-\delta}(h)(\lambda)g_\delta(V_R)^{2n}(\lambda)
\end{equation}
where $V_R(x) = V(x)\chi_{[R,\infty)}(x)$;
there is a corresponding bound for the terms arising from the multilinear
expressions $S_{2n+1}$.
We may dominate $\sup_R  g_\delta(V_R)(\lambda)^2$
by $C' g_{\delta'}(V)(\lambda)^2$ for any $\delta'>\delta$;
to simplify notation we replace $\delta'$ again by $\delta$.
Consequently by summing over $n$ and comparing with the Taylor expansion 
for the exponential function,
\begin{multline}
\Big|
\langle \int_\Lambda \int_R^\infty f(\lambda)\left( u(x,\lambda)
-e^{i \phi(x,\lambda)}\right) h(x)\,dx\,d\lambda
\rangle\Big|
\\
\le 
\int_\Lambda
Cg_{-\delta}(h)(\lambda)
g_\delta(V_R)(\lambda)
e^{Cg_\delta(V)(\lambda)^2}
|f(\lambda)|\,d\lambda.
\end{multline}

Recall that
for any $\delta>0$, $\|g_{-\delta}(h)\|_{L^2(\Lambda)}
\le C_{\Lambda,\delta}\|h\|_{L^2}\le C_{\Lambda,\delta}<\infty$,
for any compact $\Lambda\Subset(0,\infty)$, uniformly over all
martingale structures, not necessarily adapted to $h$.
The factor $\exp(Cg_\delta(V)(\lambda)^2)$ 
is bounded uniformly for $\lambda\in\Lambda$, 
by definition of a set of uniformity. 
Thus it suffices to show that 
$\|g_\delta(V_R)(\lambda)\|_{L^2(\Lambda,d\lambda)}\to 0$ as $R\to\infty$. 
Now in the sum \eqref{eq:defnofgdelta}
defining $g_\delta(V)$, the $\ell^2$ sum
over $j$ for fixed $m$ 
is $\le C2^{-\e m}\|V_R\|_{L^2}$ in $L^2(\Lambda)$ 
for some $\e>0$,
so it suffices to show that 
$\int_{E^m_j}e^{2i\phi(x,\lambda)}V_R(\lambda)\,dx
\to 0$ in $L^2(\Lambda)$. This holds by Lemma~\ref{lemma:parseval},
since $V_R\to 0$ in $L^1+L^p$ norm.
\end{proof}

\begin{corollary}  \label{cor:dropmultilinearterms}
Let $\rho>0$.
For any $f\in L^2(\reals^+,d\lambda)$ supported on $[\rho,\infty)$,
\begin{equation}
\|e^{-itH_V}U_V f-\sqrt{2/\pi}
\int_0^\infty e^{-i\lambda^2 t}f(\lambda)\psi(x,\lambda)
\,d\lambda\|_{L^2(\reals^+)}
\to 0
\text{ as } |t|\to\infty.
\end{equation}
\end{corollary}

The restriction on the support of $f$ comes about because 
of the factor of $\lambda\rp$ in the definition of $\phi$,
which causes difficulties as $\lambda\to 0$.

\begin{proof}
It is straightforward to show that 
$\int_0^\infty e^{-i\lambda^2 t}f(\lambda)\psi(x,\lambda) \,d\lambda \to 0$
in $L^2([0,R])$ for any finite $R$. 

Let $\e>0$. 
Choose a set of uniformity $\Lambda\subset[\rho,\infty)$ 
and a function $F\in L^\infty (\Lambda)$ 
such that $\|f-F\|_{L^2(\reals^+)}<\e$.
Then $\|e^{-itH_V}(U_V f-U_V F)\|_{L^2(\reals^+)}<\e$ for every
$t$, as well.
By the preceding lemma,
there exists $R<\infty$ such that
$\|e^{-itH_V}U_V F-\sqrt{2/\pi}
\int_0^\infty e^{-i\lambda^2 t}F(\lambda)\psi_0(x,\lambda)
\,d\lambda\|_{L^2([R,\infty))}<\e$ for all $t\in\reals$.
Fixing such an $R$,
there exists $T<\infty$ such that
both $e^{-itH_V}U_V F$ and $\sqrt{2/\pi}\int_0^\infty e^{-i\lambda^2 t}F(E)\psi_0(x,\lambda)
\,d\lambda$ have $L^2([0,R])$ norms $<\e$, for all $|t|\ge T$.
Thus 
\begin{equation}
\|e^{-itH_V}U_V f-\sqrt{2/\pi}\int_0^\infty e^{-i\lambda^2 t}F(\lambda)\psi_0(x,\lambda)
\,d\lambda\|_{L^2(\reals^+)}<4\e\text{ for all } |t|\ge T.
\end{equation}
Finally, note that
$\|\int_\rho^\infty g(\lambda) u_0(x,\lambda)\,d\lambda\|_{L^2(\reals^+)}
\le C_\rho \|g\|_{L^2([\rho,\infty))}$
for all $g$, where $C_\rho<\infty$ for all $\rho>0$;
this follows from the proof of Lemma~\ref{lemma:parseval}.
Applying this to $g=f-F$ allows us to replace $F$ again by $f$
in the preceding inequality, at the expense of replacing $4\e$
by $C_\rho\e$. 
Since $\e>0$ was arbitrary, this completes the proof.
\end{proof}

\section{A time-dependent phase correction} \label{section:phasereduction}

The goal of this section is to convert the leading-order asymptotics
$\exp(\pm i\lambda x \mp i(2\lambda)\rp\int_0^x V)$
to a more standard form,
$\exp(\pm i\lambda x \mp i(2\lambda)\rp\int_0^{2\lambda |t|} V)$.

\begin{lemma}  \label{lemma:easyballisticapproximation}
Let $V\in L^1+L^2$. For any $f\in L^2$ supported in a compact
subinterval of $(0,\infty)$,
\begin{multline}  \label{eq:twoeasyterms}
\Big\|\int_0^\infty e^{-i\lambda^2 t+i\lambda x 
-i(2\lambda)\rp\int_0^x V}
f(\lambda)\,d\lambda
\Big\|_{L^2(\reals^+)}
\\
+
\Big\|\int_0^\infty e^{-i\lambda^2 t+i\lambda x
-i(2\lambda)\rp\int_0^{2\lambda |t|} V}
f(\lambda)\,d\lambda
\Big\|_{L^2(\reals^+)}
\to 0 \text{ as } t\to -\infty.
\end{multline}
An analogous statement holds as $t\to+\infty$, and follows by
taking complex conjugates.
\end{lemma}

\begin{proof}
Since the mapping $f\mapsto \int_0^\infty e^{-i\lambda^2 t+i\lambda x}
 e^{-i(2\lambda)\rp\int_0^x V} f(\lambda)\,d\lambda$
is bounded from $L^2([\rho,\infty))$ to $L^2(\reals^+)$ for every
$\rho>0$, and since the variant defined by replacing
$\int_0^x V$ by $\int_0^{2\lambda t}V$ is unitary,
it suffices to prove this merely for $f\in C^\infty_0$, which we 
assume henceforth.

Let $\Lambda\Subset\reals^+$ be the support of $f$.
By integrating by parts once with respect to $\lambda$,
integrating $\exp(i[-\lambda^2 t+\lambda x])$ and differentiating
the rest, and by moving absolute value signs inside the integral, 
we obtain a pointwise in $(x,t)$ bound
\begin{multline}
|\eqref{eq:twoeasyterms}|
\le C\int_\Lambda |x-2\lambda t|\rp
\big(1+
|{\textstyle\int}_0^x V|
+ |{\textstyle\int}_0^{2\lambda |t|}V|
+ |tV(2\lambda t)| \big)
\,d\lambda
\\
\le C(|x|+|t|)\rp
(\big(1+ |{\textstyle\int}_0^x V| + |{\textstyle\int}_0^{C |t|}V| \big)
\end{multline}
for the integrands in \eqref{eq:twoeasyterms}.
For $t\le -1$,
$|x-2\lambda t|\rp\sim (x+|t|)\rp$, which tends to zero
in $L^2(\reals^+,dx)$.
Thus the first term behaves as desired. The second 
and third terms do also, if $V\in L^1$.
For $V\in L^2$, $x\rp\int_0^x V\in L^2$ as well, so 
the term $(x+|t|)\rp \int_0^x V$ is an $L^2$ function times
$x/(x+|t|)$, hence $\to 0$ in $L^2$ norm.
Lastly, $\int_0^{C|t|}V=o(|t|^{1/2})$,
while $|t|^{1/2}/(x+|t|)$ is $O(1)$ in $L^2(\reals^+)$.
\end{proof}

\begin{lemma}  \label{lemma:ballisticapproximation}
Let $V\in L^1+L^2$. For any $f\in L^2$ supported in a compact
subinterval of $(0,\infty)$,
\begin{equation}  \label{differenceofcorrections}
\Big\|
\int_0^\infty e^{-i\lambda^2 t+i\lambda x} 
\Big[ e^{-i(2\lambda)\rp\int_0^x V}-
e^{-i(2\lambda)\rp\int_0^{2\lambda |t|} V}
\Big]
f(\lambda)\,d\lambda
\Big\|_{L^2(\reals^+)}\to 0 \text{ as } t\to +\infty.
\end{equation}
\end{lemma}

\begin{proof}
As in the preceding lemma, it suffices to prove this for
$f\in C^\infty_0(\reals^+)$. 
Fix a cutoff function $\eta\in C^\infty_0(\reals)$ supported on
$(-2,2)$ and $\equiv 1$ on $[-1,1]$, with $0\le\eta\le 1$.
Let $\epsilon$ be a strictly positive function, such that
$\epsilon(t)\to 0$ as $t\to+\infty$, at a rate to be specified below. 
Set
\begin{equation}
\eta(\lambda,t,x) = \eta(\epsilon(t) t^{-1/2} (x-2\lambda t)).
\end{equation}

Consider the function of $(x,t)$
\begin{equation}  \label{eq:functionofxandt}
\int_0^\infty e^{-i\lambda^2 t+i\lambda x} 
\Big[ e^{-i(2\lambda)\rp\int_0^x V}-
e^{-i(2\lambda)\rp\int_0^{2\lambda t} V}
\Big]
f(\lambda)\,\eta(\lambda,t,x)\,d\lambda\ .
\end{equation}
We have
\begin{equation}\begin{split}
|\eqref{eq:functionofxandt}| 
&\le C \int \int_{|y-2\lambda t|\le 2\epsilon(t)\rp t^{1/2}} |V(y)|\,dy
\chi_{|x-2\lambda t|\le 2\epsilon(t)\rp t^{1/2}} \,d\lambda
\\
& \le C \int \int_{|y-x|\le 4\epsilon(t)\rp t^{1/2}} |V(y)|\,dy
\chi_{|x-2\lambda t|\le 2\epsilon(t)\rp t^{1/2}} \,d\lambda
\\
&\le C\epsilon(t)\rp t^{-1/2}
\int_{|y-x|\le 4\epsilon(t)\rp t^{1/2}} |V(y)|\,dy 
\\
&\le C\epsilon(t)^{-2}MV(x),
\end{split}\end{equation}
where $M$ is the maximal function of Hardy and Littlewood.
Observe that the restriction $\eta(\lambda,t,x)\ne 0$ for
some $\lambda\in\Lambda$ implies that $ct\le x\le Ct$
for some $c,C\in\reals^+$ depending only on $\Lambda$,
provided $\epsilon(t)\ll t^{1/2}$,
and moreover that $V$ may be replaced by its restriction to
such an interval. If $V\in L^2$ then the $L^2$ norm of
$M$ applied to the restriction of $V$ to such an interval
tends to zero as $t\to\infty$, and by choosing $\epsilon(t)$ 
to tend to zero sufficiently slowly we find that 
$\eqref{eq:functionofxandt}\to 0$ in $L^2$.
For $V\in L^1$ we have the pointwise bound
$C\epsilon(t)\rp t^{-1/2}\int_{ct}^{Ct}|V|$,
which is $o(1)\cdot\epsilon(t)\rp$  in $L^2(x\sim t)$;
this tends to zero provided $\epsilon(t)$ does so sufficiently slowly.
In the general case $V\in L^1+L^2$, we decompose 
\eqref{eq:functionofxandt}
into two parts, and estimate them separately.

There remains the contribution of $1-\eta(\lambda,t,x)$:
\begin{equation} 
\int_0^\infty e^{-i\lambda^2 t+ i\lambda x}
f(\lambda)
\Big(e^{-i(2\lambda)\rp \int_0^{x}V}
-e^{-i(2\lambda)\rp \int_0^{2\lambda t}V}
\Big)
(1-\eta)(\lambda,t,x)
\,d\lambda.
\end{equation}
We would like to apply to it the method
of stationary phase. However, the phase function
$-t\lambda^2+\lambda x-(2\lambda)\rp\int_0^{2\lambda t}V$ is not 
well behaved; its partial derivative with respect to $\lambda$ is
$x-2\lambda t
+(2\lambda)^{-2}(\int_0^{2\lambda t}V)
-(2\lambda)\rp 2tV(2\lambda t)$,
and the final term is not well under control
unless one assumes $V(x) = O(|x|^{-1/2})$.
Instead, we integrate by parts, 
integrating $\exp(i[\lambda x-\lambda^2 t])$
and differentiating the rest, to obtain
\begin{equation} \label{eq:abigjob}
i\int_0^\infty e^{-i\lambda^2 t+ i\lambda x}
f(\lambda)
\frac{\p}{\p\lambda}
\Big[(x-2\lambda t)\rp
\Big(e^{-i(2\lambda)\rp \int_0^{x}V}
-e^{-i(2\lambda)\rp \int_0^{2\lambda t}V}
\Big)
(1-\eta)(\lambda,t,x)
\Big] \,d\lambda.
\end{equation}

When the derivative is expanded according to Leibniz's rule, 
various terms result.  The main terms are those in which
$\partial/\partial\lambda$ acts on either of the two exponentials 
$\exp( i(2\lambda)\rp(\int_0^* V) )$, and we discuss these first.
One such term is a constant multiple of
\begin{multline}
\int_0^\infty e^{-i\lambda^2 t+ i\lambda x
-i(2\lambda)\rp \int_0^{2\lambda t}V}
f(\lambda)
(x-2\lambda t)\rp
(1-\eta)(\lambda,t,x)
\lambda\rp tV(2\lambda t)\,d\lambda 
\\
=
c
\int_0^\infty e^{-i(y^2/4t)+ i(xy/2t)
-it y\rp \int_0^{y}V}
\tilde f(y/2t)
(x-y)\rp
(1-\eta)(\lambda,t,x)
 v(y)
\,dy 
\end{multline}
where 
we have substituted $y = 2\lambda t$ and written
$\tilde f(\lambda) = \lambda\rp f(\lambda)\in C^\infty_0$,
and where $v=V$ (for the present moment only).
The integral operators with kernels $\exp(iAxy)(x-y)\rp 
\eta(\delta(x-y))$
are bounded on $L^2(\reals)$, uniformly in $A,\delta\in\reals^+$. 
Thus the $L^2(dx)$ norm of this last expression is $O(\|v\|_{L^2})$.
Moreover, since $f$ is supported where $\lambda\ge\rho>0$, 
only the restriction of $v$ to $[2\rho t,\infty)$ comes into play, so we 
obtain a bound of $C\|v\|_{L^2([2\rho t,\infty))}$. This holds uniformly in
all real-valued functions $V$ appearing in the exponent.

There is also an easy alternative bound $O(t^{-1/4}\|v\|_{L^1})$, 
obtained directly by inserting absolute values inside the integral.
Thus the general case $V\in L^1+L^2$ may be treated by decomposing $v$
as a sum, and estimating the two terms separately.

Another term arising from \eqref{eq:abigjob} differs 
only in that  $tV(2\lambda t)$
is replaced by $c\lambda^{-2}\int_0^{2\lambda t}V$;
in terms of the new variable $y$,
$v(y)$ is replaced by $t\rp\int_0^y v$.
Since $\lambda$ ranges over a compact interval $\Lambda$,
$y$ ranges over an interval $[ct,Ct]$ where $0<c<C<\infty$. Therefore
$t\rp\int_0^y v$ is majorized in $L^q(dy)$ by $\|v\|_{L^q(d\lambda)}$,
for every $1\le q\le\infty$.
Thus the same reasoning applies to this term.

In the last of the main terms, 
the derivative falls on $\exp(-i(2\lambda)\rp\int_0^x V)$,
and we obtain
\begin{equation}  \label{eq:slightlynastymainterm}
c\big({\textstyle\int}_0^x V\big)
\int e^{-i\lambda^2t+i\lambda x-i(2\lambda)\rp\int_0^x V}
(x-2\lambda t)\rp f(\lambda)(1-\eta)(\lambda,t,x)\,d\lambda .
\end{equation}
If $x\ge Ct$ where $C$ is a sufficiently large constant,
depending on $\Lambda$, this is majorized by
$Cx\rp\big|\int_0^x V\big|$,
which is $o(\|V\|_{L^2})$ in $L^2(|x|\ge Ct)$ as $t\to+\infty$.
However, for $x\lesssim t$ and
for $V\in L^2$, simple estimates on this quantity seem to miss the
desired conclusion by a factor of $\log(t)$, so we integrate by
parts again and move absolute values inside the integral to
obtain an upper bound
\begin{multline}
|\eqref{eq:slightlynastymainterm}|
\le C\big|{\textstyle\int}_0^x V\big|
\cdot\int_\Lambda |x-2\lambda t|^{-2}
\\
\Big(
\frac{t|1-\eta|}{|x-2\lambda t|} + |1-\eta|
+ C\epsilon t^{1/2}|\eta'(\epsilon t^{-1/2}(x-2\lambda t))|
+ \big|(1-\eta){\textstyle\int}_0^x V\big|
\Big)\,d\lambda,
\end{multline}
where $(1-\eta)\equiv (1-\eta)(\lambda,t,x)$
and $\epsilon\equiv\epsilon(t)$.
Bearing in mind that now $x\le Ct$,
one term is
\begin{multline}
\le C\big| {\textstyle\int}_0^x V\big|^2
\int_{|\lambda -(x/2t)|\ge \epsilon\rp t^{-1/2}} t^{-2}|\lambda-(x/2t)|^{-2}\,d\lambda 
\\
\le C\big| {\textstyle\int}_0^x V\big|^2
t^{-3/2}\epsilon
 \le Ct^{-3/2}\epsilon x^{1/2}\big|{\textstyle\int}_0^x V\big|
 \le C\epsilon(t) x\rp\big|{\textstyle\int}_0^x V\big|,
\end{multline}
which is $o(\|V\|_{L^2})$ in $L^2$ norm.
Another term is
\begin{equation}
\le C\big| {\textstyle\int}_0^x V\big|
\int_{|\lambda -(x/2t)|\ge \epsilon\rp t^{-1/2}} 
t^{-2}|\lambda-(x/2t)|^{-3}\,d\lambda
\le C\big| \int_0^x V\big|
\epsilon^2 t\rp,
\end{equation}
which again is $o(\|V\|_{L^2})$ in $L^2$ norm.
The two remaining terms are both majorized by the discussion of this last term,
thus completing the analysis of the case where $\partial/\partial\lambda$
acts on $\exp(-i(2\lambda)\rp\int_0^x V)$,
if $V\in L^2$.
In the general case $V\in L^1+L^2$, the contribution of the $L^1$
part is better behaved and the details are left to the reader.

Another type of term arises when $\partial/\partial\lambda$
acts in \eqref{eq:abigjob} on 
$(1-\eta)(\lambda,t,x)$, producing 
$c\epsilon(t)t^{1/2}\cdot\eta'(\epsilon(t) t^{-1/2}(x-2\lambda t)$.
The factor $\eta'(\epsilon(t) t^{-1/2}(x-2\lambda t)$
is supported where $|x-2\lambda t|\le 2\epsilon(t)\rp t^{1/2}$,
and 
$\epsilon(t)t^{1/2}/(x-2\lambda t)=O(1)$ in this region.
Therefore the analysis given above for the contribution
of $\eta(\lambda,t,x)$ applies equally well to this term.

When the derivative falls on $f(\lambda)$, we have gained 
a factor of $(x-2\lambda t)\rp$, and have an upper bound
$C\int_\Lambda |x- 2\lambda t|\rp\chi_{|x-2\lambda t|\ge 
\epsilon(t)\rp t^{1/2}}\,d\lambda$. For each $\lambda\in\Lambda$,
the integrand is $O(t^{-1/4}\epsilon(t)^{1/2})$ in $L^2(\reals^+)$.

There remains the term in which the derivative
$\partial/\partial\lambda$ falls on
$(x-2\lambda t)\rp$; this term is a constant multiple of
\begin{equation}
\int_0^\infty e^{-i\lambda^2 t+ i\lambda x}
\Big(e^{-i(2\lambda)\rp \int_0^{2\lambda t}V}
-e^{-i(2\lambda)\rp \int_0^{x}V}
\Big)
f(\lambda)
\frac{t}{(x-2\lambda t)^{2}}
(1-\eta)(\lambda,t,x)
\,d\lambda
\end{equation}
and hence is majorized by
$C\int_{\Lambda} t|x-2\lambda t|^{-2}\chi_{|\lambda-(x/2t)|\ge\epsilon(t)\rp t^{-1/2}}
\,d\lambda\le C\epsilon(t) t^{-1/2}$,
which gives the desired $L^2$ bound for $x\le C_0t$,
for any fixed $C_0$.
Choosing $C_0$ to be sufficiently large,
we can more simply majorize $|x-2\lambda t|^{-2}$ 
for $x\ge C_0 t$ by $Cx^{-2}$ in the integrand,
obtaining the pointwise bound $C t/x^2$ for the integral,
giving a contribution $\le C t^{-1/2}$ to the $L^2(dx)$ norm.
\end{proof}

\section{Modified wave operators} \label{section:modifiedwave}

\subsection{The half-line case}

We discuss modified wave operators for Schr\"odinger operators
on the half line. 

Recall from the introduction the modified wave operators 
\eqref{modifiedwaveoperatorsdefn}
\[ \Omega^m_\pm f
= \lim_{t\to\mp\infty} e^{itH_V}e^{-it H_0 \pm iW(H_0,\mp t)}f \]
for all $f\in L^2(\reals^+)$, where existence of the limit remains
to be established. Here $W$ is given by \eqref{phasec}: 
\[ W(\lambda,t) = -(2\lambda)\rp\int_0^{2\lambda t} V(s)\,dx. \]
Motivating this definition is the heuristic that a wave packet
with frequency $\approx\lambda$ should propagate
with velocity $\approx 2\lambda$,
so that for large $t$, $\int_0^x V$ should behave like
$\int_0^{2\lambda t}V$.

\begin{proof}[Proof of Theorem~\ref{thm:waveoperators}]
We write $A(t)\sim B(t)$ to mean that
$\|A(t)-B(t)\|_{L^2(\reals^+)}\to 0$ as $t\to+\infty$,
or as $t\to -\infty$, depending on which limit is presently
under discussion.
Given $g\in L^2(\reals^+)$, write 
\begin{equation}
g(x)
= \pi\rp
\int_0^\infty \hat g(\lambda) 
[e^{i\lambda x}-e^{-i\lambda x}]
\,d\lambda,
\end{equation}
where this defines the modified Fourier transform
$\hat g$, so that  $\|\hat g\|_{L^2(\reals^+)}
= \sqrt{\pi}\|g\|_{L^2(\reals^+)}$.
Then
\begin{equation}
e^{-itH_0}g(x)  
= \pi\rp\int_0^\infty \hat g(\lambda)e^{-i\lambda^2t}
[e^{i\lambda x}-e^{-i\lambda x}]
\,d\lambda,
\end{equation}
and the modified evolution is
\begin{equation}
e^{-itH_0 -iW(H_0,t)}g(x)
= \pi\rp
\int_0^\infty \hat g(\lambda) e^{-i\lambda^2 t}
\Big[e^{i\lambda x-i(2\lambda)\rp\int_0^{2\lambda t}V}
-e^{-i\lambda x-i(2\lambda)\rp\int_0^{2\lambda t}V}\Big]
\,d\lambda.
\end{equation}

Define $\tilde g\in\scripth_{\text{ac}}$ by
\begin{equation}
\tilde g(x)=\pi\rp\int_0^\infty \hat g(\lambda) \psi(x,\lambda)\,d\lambda.
\end{equation}
The map $g\mapsto \tilde g$ is a unitary bijection from
$L^2(\reals^+)$ to $\scripth_{\text{ac}}$.
We aim to prove that $e^{-itH_0 -W(H_0,t)}g
\sim e^{-itH_V}\tilde g$, as $t\to+\infty$, for any $g\in L^2(\reals^+)$. 
By unitarity of the evolutions, this implies the existence of 
the modified wave operator $\Omega^m_-$, and likewise that
it is an isometric bijection from
$L^2(\reals^+)$ to $\scripth_{\text{ac}}$.
Moreover, by unitarity, it suffices to prove this 
under the assumption that $\hat g\in C^\infty_0(\reals^+)$,
which we assume henceforth.

Now
\begin{align*}
e^{-itH_0 -iW(H_0,t)}g(x)
&=
\pi\rp\int_0^\infty
\hat g(\lambda)e^{-i\lambda^2 t}
\Big[e^{i\lambda x-i(2\lambda)\rp\int_0^{2\lambda t}V}
-
e^{-i\lambda x-i(2\lambda)\rp\int_0^{2\lambda t}V}
\Big]
\,d\lambda
\\
& \sim
\pi\rp\int_0^\infty
\hat g(\lambda)e^{-i\lambda^2 t}
e^{i\lambda x-i(2\lambda)\rp\int_0^{2\lambda t}V}
\,d\lambda
\end{align*}
by taking complex conjugates in
Lemma~\ref{lemma:easyballisticapproximation}.
On the other hand,
\begin{align*}
e^{-itH_V}\tilde g
&=
\pi\rp\int_0^\infty
\hat g(\lambda) e^{-i\lambda^2 t}
\psi(x,\lambda)\,d\lambda
\\
&=
\pi\rp\int_0^\infty
\hat g(\lambda) e^{-i\lambda^2 t}
\big(u(x,\lambda^2)- \frac{\bar\gamma(\lambda^2)}{\gamma(\lambda^2)}
\ov{u(x,\lambda^2)}\big)
\,d\lambda
\\
&\sim
\pi\rp\int_0^\infty
\hat g(\lambda) e^{-i\lambda^2 t}
\big(e^{i\lambda x-i(2\lambda)\rp\int_0^x V}
-
\frac{\bar\gamma(\lambda^2)}{\gamma(\lambda^2)}
e^{-i\lambda x+i(2\lambda)\rp\int_0^x V}\big)
\,d\lambda
\\
&\sim
\pi\rp\int_0^\infty
\hat g(\lambda) e^{-i\lambda^2 t}
e^{i\lambda x-i(2\lambda)\rp\int_0^x V}
\,d\lambda.
\\
&\sim
\pi\rp\int_0^\infty
\hat g(\lambda) e^{-i\lambda^2 t}
e^{i\lambda x-i(2\lambda)\rp\int_0^{2\lambda t} V}
\,d\lambda.
\end{align*}
by first Corollary~\ref{cor:dropmultilinearterms},
then Lemma~\ref{lemma:easyballisticapproximation},
and then Lemma~\ref{lemma:ballisticapproximation}.
Thus we have the asymptotic relation 
$e^{-itH_0 -iW(H_0,t)}g \sim e^{-itH_V}\tilde g$, as $t\to+\infty$, 
so $\Omega^m_- = U_V^{-1}\circ U_0$.

The analysis as $t\to-\infty$ is the same, except for the appearance
of the unimodular factor $\ov{\gamma(\lambda^2)}/\gamma(\lambda^2)$.
The conclusion is therefore that
$\Omega^m_+ = U_V^{-1}\circ[\ov{\gamma(\lambda^2)}/\gamma(\lambda^2)]
\circ U_0$, where the inner factor denotes the operator defined
on $L^2(\reals^+,d\lambda)$
by multiplication by this function of $\lambda$.

Composing, we find that $S^m = (\Omega^m_-)\rp\Omega^m_+
= U_0^{-1}\circ[\ov{\gamma(\lambda^2)}/\gamma(\lambda^2)]\circ U_0$.
\end{proof}

\begin{proof}[Proof of Theorem~\ref{mollerwo}]
Suppose now that the improper integral $\int_0^\infty V(y)\,dy$
exists. The asymptotic behavior of $V$ enters both into the
definition of the modified evolution $e^{-itH_0 -W(H_0,t)}$,
and into the definition of the scattering coefficients $\gamma(E)$.
To sort this out, we change the definition of the phase
$\xi(x,\lambda^2)$ to $\lambda x$.
We correspondingly  modify the normalization of the 
generalized eigenfunctions $u(x,\lambda^2)$ to $u(x,\lambda^2)\sim e^{i\lambda x}$. 
The scattering coefficient $\gamma(\lambda^2)$ is now defined,
for almost every $\lambda>0$,
by the relation \eqref{reflectioncoefficientdefn}
$\gamma(\lambda^2) = 1/\psi(0,\lambda)$.
Writing $\gamma(\lambda^2)= |\gamma(\lambda^2)|e^{i\arg{\gamma(\lambda^2)}}$,
we change the definition of $\psi(x,\lambda)$ to 
\begin{equation}
\psi(x,\lambda) = e^{i\arg{\gamma(\lambda^2)}}  u(x,\lambda^2)
- e^{-i\arg{\gamma(\lambda^2)}}\ov{u(x,\lambda^2)}.
\end{equation}
In the formal expressions for the spectral projectors and wave group,
$\psi$ is now replaced by this new $\psi$. 
Retracing the above analysis, we find that the wave operators exist,
and
\begin{equation}\begin{split}
\Omega_- &= U_V^{-1}\circ e^{-i\arg{\gamma(\lambda^2)}}
\circ U_0,
\\
\Omega_+ &= U_V^{-1}\circ e^{i\arg{\gamma(\lambda^2)}}
\circ U_0, \label{mwo}
\\
S &= U_0\rp\circ e^{2i\arg{\gamma(\lambda^2)}}
\circ U_0.
\end{split}\end{equation}
\end{proof}

{\em Remark.\/}
The exact expressions for the modified wave operators $\Omega^m_\pm$
depend on a choice of normalization of the solutions $\psi$; see below.
These solutions can be modified by factors $e^{i\kappa(\lambda)}$, leading
to different $U_V$ and hence different looking expressions for the wave
operators, like in \eqref{mwo} and in the proof of Theorem~\ref{thm:waveoperators}. 
However, the scattering matrix $S(\lambda)$ is invariant
under choice of such normalization.

\subsection{The whole-line case}

In the full-axis Schr\"odinger case, the results are
similar. One defines modified wave operators by
\begin{equation}  \label{eq:wholeaxiswaveoperatorsdef}
\Omega_\pm^m f =\lim_{ t \rightarrow \mp \infty}
e^{itH_V}e^{-iW_a(-i\partial_x,t)}  f, 
\end{equation}
where $\partial_x$ is the operator of one differentiation in $x$, and
\begin{equation} W_a(\lambda,t) = \lambda^2 + \frac{1}{2\lambda} 
\int_0^{2\lambda t} V(s)\,ds. \end{equation}
The modified free evolution operator can be written as 
(compare with \eqref{frevwa})
\begin{equation}\label{modevol}
e^{-iW_a(-i\partial_x,t)}f(x) =
\frac{1}{2\pi} \int_{-\infty}^{\infty} e^{-i\lambda ^2t +i\lambda x -
\frac{i}{2\lambda} \int_0^{2\lambda t} V(s)\,ds} \hat{f}(\lambda)\,d\lambda,
\end{equation}
where $\hat{f}(\lambda)= \int \exp(-i\lambda x)f(x)\,dx$ 
is the Fourier transform of $f$.
Recall the representation \eqref{perevwa} for the perturbed evolution
(we set here $E=\lambda^2$)
\begin{equation}\label{impevol}
 e^{-iH_Vt} f(x) = \frac{1}{2\pi} \int_0^\infty e^{-i\lambda^2t}
(\psi_+(x,\lambda)\tilde{f}_+(\lambda) 
+ \psi_-(x,\lambda)\tilde{f}_-(\lambda) )\,d\lambda.
\end{equation}
Denote by  $U_V$ the operator $f\mapsto (f_+,f_-)$, where
$f_\pm(\lambda) = \int_\reals \ov{\psi}_\pm (x,\lambda) f(x)\,dx$.
Denote by $\scripth_{\text{ac}}$ the maximal closed subspace
of $\scripth=L^2(\reals)$ on which $H_V$ has purely absolutely
continuous spectrum.

\begin{theorem}
Let $V\in L^1+L^p(\reals)$ for some $1<p<2$.
Then for every $f\in L^2$, the limits in 
\eqref{eq:wholeaxiswaveoperatorsdef}
exist in $L^2(\reals)$ norm as $t\to\mp\infty$.
The modified wave operators $\Omega^m_\pm$ thus defined are 
surjective and unitary from $L^2(\reals)$ to $\scripth_{\text{ac}}$.
One has 
\begin{equation}\begin{split}
\Omega^m_+ &=U_V^{-1} U_0,
\\
\Omega^m_- &= U_V^{-1} S(\lambda)^{-1} U_0,
\\
S^m &= U_0\rp S(\lambda) U_0
\end{split} \end{equation}
where $S(\lambda)$ denotes multiplication by the scattering matrix
\begin{equation}\label{scatmat}
S(\lambda)
= \left( \begin{array}{cc} t_1(\lambda)
& - \ov{r}_1(\lambda) \frac{t_1(\lambda)}{\ov{t}_1(\lambda)} \\
r_1(\lambda) & t_1(\lambda) \end{array} \right).
\end{equation}
\end{theorem}

As in the half-line case, we could make the wave operators look more
symmetric by modifying the definition of $U_V$. Let $A = \sqrt{S}$
be any matrix square root of the unitary operator $S$. Then if we were
to define $\tilde U_V = AU_V$,  we would obtain
$\Omega^m_+ = \tilde U_V\rp A^* U_0$,
and
$\Omega^m_- = \tilde U_V\rp A U_0$.


We sketch the asymptotic analysis of
$e^{-itH_V}g$, as $|t|\to\infty$, 
for arbitrary $g\in\scripth_{\text{ac}}$.
We proceed formally; all steps are justified as in the preceding
subsection.
Write
$g = (2\pi)\rp \int_0^\infty
\big[
\hat g_+(\lambda)\psi_+(x,\lambda)
+ \hat g_-(\lambda)\psi_-(x,\lambda)
\big]
e^{-i\lambda^2 t}\,d\lambda$
with $2\pi \|g\|_{L^2}^2 = \|\hat g_+\|_{L^2(\reals^+)}^2
+ |\hat g_-\|_{L^2(\reals^+)}^2$.
Define 
\begin{equation}\begin{split}
\Phi_+(\lambda,x,t) &=
\phantom{-}\lambda x -(2\lambda)\rp\textstyle\int_0^{2\lambda t}V,
\\
\Phi_-(\lambda,x,t) &=
-\lambda x +(2\lambda)\rp\textstyle\int_0^{-2\lambda t}V .
\end{split}\end{equation}
Formally,
\begin{equation}
e^{-itH_V}\psi_+(x,t)
\sim 
\left\{
\begin{aligned}
\big[
&t_1(\lambda)e^{i\Phi_+(\lambda,x,t)}
+ r_1(\lambda)e^{i\Phi_-(\lambda,x,t)}
\big]
e^{-i\lambda^2 t}
&\ \text{ as } t\to+\infty,
\\
&e^{i\Phi_+(\lambda,x,t)}e^{-i\lambda^2 t}
&\ \text{ as } t\to-\infty,
\end{aligned}
\right.
\end{equation}
and
\begin{equation}
e^{-itH_V}\psi_-(x,t)
\sim 
\left\{
\begin{aligned}
\big[
&t_2(\lambda)e^{i\Phi_-(\lambda,x,t)}
+ r_2(\lambda)e^{i\Phi_+(\lambda,x,t)}
\big]
e^{-i\lambda^2 t}
&\ \text{ as } t\to+\infty,
\\
&e^{i\Phi_-(\lambda,x,t)}e^{-i\lambda^2 t}
&\ \text{ as } t\to-\infty,
\end{aligned}
\right.
\end{equation}
Therefore as $t\to+\infty$,
\begin{multline}
e^{-itH_V}(2\pi)\rp\int_0^\infty
\big[
\hat g_+(\lambda)\psi_+(x,\lambda)
+ \hat g_-(\lambda)\psi_-(x,\lambda)
\big]
e^{-i\lambda^2 t}\,d\lambda
\\
\sim
(2\pi)\rp\int_0^\infty
\Big(
\hat g_+(\lambda)
\big[
t_1(\lambda)e^{i\Phi_+(\lambda,x,t)}
+ r_1(\lambda)e^{i\Phi_-(\lambda,x,t)}
\big]
\\
+ \hat g_-(\lambda)
\big[
t_2(\lambda)e^{i\Phi_-(\lambda,x,t)}
+ r_2(\lambda)e^{i\Phi_+(\lambda,x,t)}
\big]
\Big)
e^{-i\lambda^2 t}\,d\lambda .
\end{multline}
As $t\to-\infty$, we have instead the asymptotics
\begin{equation}
\sim
(2\pi)\rp\int_0^\infty
\big(
\hat g_+(\lambda)
e^{i\Phi_+(\lambda,x,t)}
+ \hat g_-(\lambda)
e^{i\Phi_-(\lambda,x,t)}
\big)
e^{-i\lambda^2 t}\,d\lambda .
\end{equation}
Here ``$\sim$'' signifies asymptotic equality in $L^2(\reals,dx)$.
Thus for any function $f(\lambda) = (f_+(\lambda),f_-(\lambda))^T$
taking values in $L^2(\reals^+,\complex^2)$,
\begin{equation}
e^{-itH_V}U_V\rp \begin{pmatrix}f_+ \\ f_- \end{pmatrix}
\sim
\left\{\begin{aligned}
e^{-iW_a(-i\partial_x,t)}U_0\rp
\begin{pmatrix}
t_1(\lambda) & r_2(\lambda)
\\ r_1(\lambda) & t_2(\lambda)
\end{pmatrix}
\begin{pmatrix}f_+ \\ f_- \end{pmatrix}
&\ \text{ as } t\to+\infty, \\
e^{-iW_a(-i\partial_x,t)}U_0\rp
\begin{pmatrix}
1 & 0 \\ 0 & 1
\end{pmatrix}
\begin{pmatrix}f_+ \\ f_- \end{pmatrix}
&\ \text{ as } t\to-\infty.
\end{aligned}\right.
\end{equation}
Inverting this and exploiting the identities relating the various
scattering coefficients to simplify the resulting formula, we obtain
\begin{equation}
e^{-iW_a(-i\partial_x,t)}U_0\rp
\begin{pmatrix}f_+ \\ f_- \end{pmatrix}
\sim
\left\{\begin{aligned}
e^{-itH_V}U_V\rp 
\begin{pmatrix}
\bar t_1 & \bar r_1
\\ -(r_1\bar t_1/t_1) & \bar t_1
\end{pmatrix}
\begin{pmatrix}f_+ \\ f_- \end{pmatrix}
&\ \text{ as } t\to+\infty,
\\
e^{-itH_V}U_V\rp 
\begin{pmatrix}
1 & 0 \\ 0 & 1
\end{pmatrix}
\begin{pmatrix}f_+ \\ f_- \end{pmatrix}
&\ \text{ as } t\to-\infty.
\end{aligned}\right.
\end{equation}
with asymptotic equality holding in $L^2(\reals^+,\complex^2,d\lambda)$
norm.
The conclusions of the theorem follow directly.

As in the half-line case, we obtain a cleaner conclusion under 
the additional hypothesis that the improper integrals
$\int_0^{\pm\infty}V(x)\,dx$ exist.
Modify the definitions of the scattering coefficients by redefining
\begin{align}
&\psi_+(x,\lambda) = \left\{ \begin{array}{ll} t_1(\lambda) e^{i\lambda x} \ + o(1),
& x \rightarrow +\infty, \\
\left( e^{i \lambda x} +r_1(\lambda)e^{-i \lambda x} \right)\ + o(1),
& x \rightarrow -\infty \end{array} \right.
\label{risc2}
\\
& \psi_-(x,\lambda) = \left\{ \begin{array}{ll} t_2(\lambda) e^{-i\lambda x} \ + o(1),
& x \rightarrow -\infty, \\
\left( e^{-i\lambda x} +r_2(\lambda)e^{i \lambda x} \right)\ + o(1),
& x \rightarrow +\infty. \end{array} \right.
\label{lesc2}
\end{align}

\begin{theorem}
Let $V\in L^1+L^p(\reals)$ for some $1<p<2$, 
and suppose that both the improper integrals
$\int_0^{\pm\infty}V(x)\,dx$ exist. Then for every $f\in L^2(\reals)$,
the two limits
\begin{equation}
\Omega_\pm f = \lim_{t\to\mp\infty}
e^{itH_V}e^{-itH_0}f
\end{equation}
exist in $L^2(\reals)$ norm, and the operators thus defined are
surjective and unitary from $L^2(\reals)$ to $\scripth_{\text{ac}}$.
One has 
\begin{equation}\begin{split}
\Omega_+ &=U_V^{-1} U_0,
\\
\Omega_- &= U_V^{-1} S(\lambda) U_0,
\\
S &= U_0\rp S(\lambda) U_0
\end{split} \end{equation}
where $S(\lambda)$ is defined in terms of the modified
scattering coefficients $t_j(\lambda),r_j(\lambda)$ defined 
in \eqref{risc2},\eqref{lesc2}, and $U_V$  is defined as before 
but in terms of $\psi_\pm$ as in \eqref{risc2}, \eqref{lesc2}.
\end{theorem}

\subsection{The Dirac case}  \label{subsection:dirac}

Consider next the operator $D_V$ on $\scripth=L^2(\reals)$.
Because no phase correction $\int_0^x V$ appears in the generalized
eigenfunction asymptotics, there is no need for modification of the standard
wave operators in this case.

We calculate 
\begin{multline}
e^{-itD_V}U_V\rp G
\sim 4\rp\int_\reals
e^{-iEt}
\Big[
G_+(E)\Big(
t_1(E)e^{iEx}\begin{pmatrix} 1 \\ i\end{pmatrix}
+ r_1(E)e^{-iEx}\begin{pmatrix} 1 \\ -i\end{pmatrix}\Big)
\\
+
G_-(E)\Big(
t_2(E)e^{-iEx}\begin{pmatrix} 1 \\ -i\end{pmatrix}
+ r_2(E)e^{-iEx}\begin{pmatrix} 1 \\ i\end{pmatrix}\Big)
\Big]
\,dE
\end{multline}
as $t\to+\infty$,
and
\begin{equation}
\sim 4\rp\int_\reals
e^{-iEt}
\Big(
G_+(E)e^{iEx}\begin{pmatrix} 1 \\ i\end{pmatrix}
+ G_-(E)e^{-iEx}\begin{pmatrix} 1 \\ -i\end{pmatrix}
\Big)
\,dE
\end{equation}
as $t\to-\infty$.

\begin{theorem}
Let $V\in L^1+L^p(\reals)$ for some $1<p<2$.
Then for each $f\in L^2(\reals)$, both of the limits
\begin{equation}
\Omega_\pm f = \lim_{t\to\mp\infty} e^{-itD_V}e^{itD_0}f
\end{equation}
exist in $L^2(\reals)$ norm, and define surjective and unitary
operators to $\scripth_{\text{ac}}$.
Moreover $\Omega_+ = U_V^{-1} U_0$ and $\Omega_- = U_V^{-1} S(E)^{-1}U_0,$ 
where $S(E)$ is defined as before by \eqref{scatmat}. 
The scattering operator $S$ is equal to 
$U_0^{-1} S(E) U_0.$
\end{theorem}

\section{Asymptotic completeness}  \label{section:completeness}

Let us denote by $\scripth_{\text{pp}}(A)$ the pure point subspace of the self-adjoint operator $A.$
 Recall from \cite{RS3}
\begin{definition}
The (modified) wave operators $\Omega_\pm$ are called asymptotically complete if their 
ranges coincide with the  orthogonal complement of $\scripth_{\text{pp}}(H_V)$. 
\end{definition}
In our context, since the absolutely continuous spectrum is under control, we have 
\begin{corollary}\label{sccor}
The wave operators of Theorems~\ref{thm:waveoperators} and \ref{wo} are asymptotically 
complete if and only if the singular continuous spectrum of the operator $H_V$ is empty.
\end{corollary}

In general, the question of asymptotic completeness of the wave operators for potentials in 
$L^p$ remains open for $1<p<2$. 
There exist examples with embedded dense point spectrum on $\reals^+$ 
\cite{Nab,simondenseembedded}.
A very recent preprint \cite{denisov} states that singular continuous spectrum can arise
for $L^2$ potentials.
There also exist estimates on the size of the set where the singular spectrum may be 
supported \cite{Re2,christkiselevslowlyvarying} and  examples of operators which have decaying solutions on a half axis for a set of spectral parameters having 
exactly the right dimension \cite{Re3,KR}. However,
no examples have yet been constructed of Schr\"odinger operators, with 
potential $V \in L^p$ for some $p< 2$, possessing nonempty singular continuous spectrum. 

Nevertheless, there are some settings where generic (in some sense) asymptotic completeness can be 
inferred.  We discuss two such cases: almost surely for certain random models, 
and for almost every boundary condition with any fixed potential.
In both of these cases, the spectrum is purely absolutely continuous where $E>0$,
so asymptotic completeness in the relatively weak sense in which we have defined it
implies a stronger form.

Let us say that the half-axis operator $H_V^\beta$ has boundary condition $\beta$ at the origin if 
$u(0)+\beta u'(0)=0$ for any function $u$ in the domain of $H_V^\beta.$   
\begin{theorem}\label{bcascom}
Let $H_V^\beta$ be a Schr\"odinger operator defined on a half-axis with the boundary condition $\beta.$ 
Assume that $V \in L^1+L^p$ for some $1<p<2$. 
Then for almost every $\beta$, the modified wave operators defined in 
\eqref{modifiedwaveoperatorsdefn} are asymptotically complete.  
Moreover if $\int_0^\infty V(x)\,dx$ exists, then the usual wave operators 
defined by \eqref{wo} are asymptotically complete. 
\end{theorem}
\begin{proof}
The result is a simple corollary of Theorems~\ref{thm:waveoperators}, \ref{mollerwo},
and general rank one perturbation theory. 
Theorem~\ref{thm:waveoperators} and well-known results of scattering theory imply that the absolutely continuous 
part of the spectral measure  of $H_V$ corresponding to some boundary condition $\beta_0,$ 
$\mu_{\rm{ac}}^{\beta_0},$ fills all of $R^{+};$ 
that is, $D\mu_{\rm{ac}}^{\beta_0}(E)>0$ for a.e. $E.$ The variation of the boundary condition can 
be regarded as rank one perturbation \cite{SVan}, Section I.6. The standard rank one perturbation theory 
then implies that for any $\beta$, the singular 
part of the spectral measure $\mu_{\rm{s}}^{\beta}$
can only be supported on a fixed set of energies $S \subset \reals^+$ of zero Lebesgue measure 
\cite{SVan}, Theorem II.2. But then again by rank one theory for almost every
$\beta$ the singular spectrum on $\reals^+$ is empty \cite{SVan}, Theorem I.8. 
Since the potential belongs to $L^1+L^p,$ the spectrum below zero is 
discrete (with only $0$ as a possible accumulation point).   
\end{proof}
Next, let us consider the following random model:
\begin{equation}\label{rex}
V(x) = \sum\limits_{n=1}^\infty a_n(\omega)g(n) f(x-n),
\end{equation}
where $f(x) \in C_0^\infty (0,1),$ $g \in l^p$ for some $p<2$, 
and $a_n(\omega)$ are independent identically 
distributed bounded random variables with zero expectation.  We have
\begin{theorem}\label{ranascom}
Let $H_V$ be a one dimensional Sch\"odinger operator with random potential  \eqref{rex}. Then 
with probability one, the wave operators $\Omega_\pm$ exist and are asymptotically complete.  
\end{theorem}
Indeed, Theorem 9.1 of \cite{KLS} shows that almost surely, 
the spectrum of the Schr\"odinger operator with 
potential \eqref{rex} is purely absolutely continuous on $\reals^+.$ 
Notice that our assumptions on the potential easily imply
that the improper integral $\int_0^\infty V$ exists almost surely.
Theorem~\ref{thm:waveoperators}
and decay of the potential then imply asymptotic completeness. 

This illustrates another type of situation 
where asymptotic completeness holds. We remark that the result holds in a variety of more general 
random models for which \eqref{rex} is just an illustration. 
For a more general setting, see \cite{KLS}. 

\smallskip
{\it Discussion.\/}
The above is only one of various possible definitions of asymptotic completeness. 
The notion of asymptotic 
completeness is intended to describe a situation where the Hilbert space is split into 
two orthogonal subspaces \cite{RS3}: $\scripth_{\text{pp}}(H_V)$ and the range of wave operators, 
$\scripth_{\text{ac}}(H_V).$ 
On $\scripth_{\text{ac}}(H_V)$ the perturbed dynamics is close to the modified free evolution 
at large times, and corresponds to the scattering states. On  $\scripth_{\text{pp}}(H_V)$ the dynamics
is supposed to be bounded in some sense. However, the intuitive physical assumption 
that pure point 
spectrum leads to dynamics which is bounded needs to be clarified, 
and in recent years there have appeared examples with very non-trivial transport on the 
pure point component. A 
widely accepted way to calibrate transport properties is to consider evolution of 
the averaged moments of coordinate operator: 
\begin{equation}\label{moments}
\langle \langle |X|^m \rangle_\phi \rangle_T = \frac{1}{T} \int\limits_0^T | \langle e^{-iHt}\phi, |X|^m e^{-iHt}\phi
\rangle | \, dt, 
\end{equation}
 where $\phi$ is the initial state, $\langle \phi_1, \phi_2 \rangle$ is the inner product and $|X|^m$
the operator of multiplication by $(|x|+1)^m$ in coordinate representation. 
The paper \cite{DJLS} contains an example of a (discrete) Schr\"odinger operator $h$ with pure point spectrum and 
exponentially decaying eigenfunctions, such that 
\begin{equation}\label{pointran}
 \limsup_{t \rightarrow \infty} 
\langle \langle |X|^2 \rangle_{\delta_0} \rangle_T /T^{\alpha} = \infty 
\end{equation}
for any $\alpha <2.$ Here the initial state is $\delta_0,$  the vector localized at the origin. Given that the rate
of growth $T^2$ for the second moment corresponds to ballistic motion, as for the free Laplace
operator, this example shows that in some sense 
the transport associated with point spectrum can be 
very fast.

However, there is still an important difference between transport 
associated with point and singular continuous spectrum. Namely, let $B_R$ denote the ball of radius $R$
centered at the origin and $B_R^c$ its complement. Let $P_{{\rm pp}}$ and $P_{{\rm c}}$ be 
the orthogonal projections 
on $\scripth_{\text{pp}}(H_V)$ and the continuous subspace $\scripth_{\text{c}}(H_V)$
respectively.
Then 
for any $\epsilon >0$ there exists $R_\epsilon$ such that 
\begin{equation}\label{ppbound}
 \|e^{-iH_Vt} P_{{\rm pp}}\phi \|^2_{L^2(B^c_{R_\epsilon})} < \epsilon 
\end{equation}
for all $t.$ The growth of the moments in \eqref{pointran}
is achieved not because of the motion 
of the whole wavepacket, but because of thin tails escaping to infinity. 
On the other hand, 
we have
\begin{equation}\label{probfin}
 \frac{1}{T} \int\limits_0^T \| e^{-iH_Vt} P_{{\rm c}}
\phi \|^2_{L^2(B_R)} \,dt \stackrel{T \rightarrow \infty}{\longrightarrow} 0
\end{equation}
for any finite $R.$ 
Equation \eqref{probfin} is one of the statements of
the RAGE theorem (see, e.g. \cite{CFKS}) and is
basically a corollary of Wiener's theorem 
on Fourier transforms of measures. 
Morever, there exist examples of Schr\"odinger operators \cite{KL} in which 
the dynamics corresponding to the singular continuous subspace is almost ballistic
in a sense that the whole wavepacket is moving to infinity at a fast rate: 
for any $\rho>0$ there exists $C_\rho$ such that 
\[ \frac{1}{T} \int\limits_0^T \| e^{-iH_Vt}
\phi \|^2_{L^2(B_{C_\rho T^{1-\epsilon}})} \, dt < \rho \]
for all $T$ and $\phi$ lying in the singular continuous subspace of $H_V$. 

\section{Potentials in $\ell^p(L^1)$}

Following \cite{christkiselevdecaying}, we sketch here the small modifications needed
to extend the analysis of potentials in $L^p$ to those in $L^1+L^p$,
and indeed those in the larger class $\ell^p(L^1)$.
A locally integrable function $f$ is said to belong to the amalgamated space
$\ell^p(L^1)$ if 
\begin{equation*}
\sum_{n\in\integers} \big(\int_n^{n+1}|f(x)|dx \big)^p<\infty.
\end{equation*}
The norm $\|f\|_{\ell^p(L^1)}$ is the $p$-th root of this expression.
For any $1\le p<\infty$, this defines the
Banach space $\ell^p(L^1)$, which contains $L^1+L^p$.

A martingale structure $\{E^m_j\}$ is said to be adapted to $f$ in
$\ell^p(L^1)$ if 
\begin{equation}
\|f\cdot\chi_{E^m_j}\|_{\ell^p(L^1)}^p 
\le 2^{-m} \|f\|_{\ell^p(L^1)}^p
\end{equation}
for all $m,j$.
For any $f \in\ell^p(L^1)$, there does exist an adapted martingale
structure \cite{christkiselevdecaying}.

Lemma~\ref{lemma:parseval} extends to $\ell^p(L^1)$: this Banach space is mapped boundedly
to $L^{p'}(\Lambda)$ for any compact interval $\Lambda\Subset(0,\infty)$,
for any $1\le p\le 2$. The proof is essentially unchanged; see  
the analogous proof of Proposition~3.5 of \cite{christkiselevdecaying}.

Corollary~\ref{cor:nontangentialmaxbound} may now be refined by replacing $\|f\|_{L^p}$
by $\|f\|_{\ell^p(L^1)}$ on the right-hand side of each inequality.
With these bounds for the operator $G$ in hand, the remainder of the proof is unchanged.

\medskip
{\bf Acknowledgement} We thank Barry Simon for useful discussions. \\

\end{document}